# The Riemann Hypothesis for the Goss zeta function for $\mathbb{F}_q[T]$

Jeffrey T. Sheats

Department of Mathematics, University of Arizona,

Tucson, Arizona 85721

December 30, 1997

**Abstract**

Let $q$ be a power of a prime $p$. We prove an assertion of Carlitz which takes $q$ as a parameter. Diaz-Vargas' proof of the Riemann Hypothesis for the Goss zeta function for $\mathbb{F}_p[T]$ depends on his verification of Carlitz's assertion for the specific case $q = p$ [D-V]. Our proof of the general case allows us to extend Diaz-Vargas' proof to $\mathbb{F}_q[T]$.

In [Gos1] Goss presents zeta functions for function fields of finite characteristic which possess many interesting properties including analogs of properties of the Riemann zeta function. In this paper we prove an analog of the Riemann Hypothesis for the Goss zeta function for $\mathbb{F}_q[T]$ where $q$ is a power of a prime $p$. The statement of the analog appears in Theorem 1.1 below. For the case $q = p$ the analog was proven by D. Wan [Wan]. Our proof follows another proof for the case $q = p$ recently put forward by Diaz-Vargas [D-V]. A key part of Diaz-Vargas' proof involves a result of Carlitz concerning the vanishing of certain power sums [Car3]. The proof that Carlitz sketched for his result depends on a combinatorial assertion which takes $q$ as a parameter. Carlitz gives no justification for this assertion. A restatement of Carlitz's assertion appears as our Theorem 1.2. Diaz-Vargas proved the assertion for $q = p$. This enabled him to construct his elegant proof of the Riemann Hypothesis for $\mathbb{F}_p[T]$. Poonen proved Carlitz's assertion for the case $q = 4$ [Po]. Our main result is a proof of Carlitz's assertion for all $q$.

This proof was inspired by [Po] and uses some of the same ideas. Carlitz uses his assertion in the proofs of two results. As they are now fully justified, we restate them in Theorem 1.4.

In Section 1 we define the Goss zeta function for $\mathbb{F}_q[T]$ and then state Theorem 1.1: the analog of the Riemann Hypothesis. Next we state Carlitz's result and Theorem 1.2: the combinatorial assertion on which the result depends. A brief outline of the proof of Theorem 1.2 is given at the end of Section 1. In Section 2 we prove Theorem 1.1. Sections 3 through 7 are devoted to the proof of Theorem 1.2.

The author (a combinatorialist) wishes to thank Dinesh Thakur and David Goss for their help with the number theory. The author thanks Dinesh Thakur for bringing Carlitz's assertion to his attention. Special thanks go to the referee who made valuable stylistic suggestions and also found many mistakes in earlier versions of this paper.

## 1. Main Results

In this paper $\mathbb{N}$ denotes the set of nonnegative integers and we set $\mathbb{Z}^+ := \mathbb{N}\setminus\{0\}$. Let $p$ be prime and set $q := p^s$. Let $v$ denote the $T^{-1}$-adic valuation on $K := \mathbb{F}_q(T)$. Then the field of Laurent series $K_\infty := \mathbb{F}_q((T^{-1}))$ is the completion of $K$ with respect to $v$. Denote by $A^+$ the set of monic polynomials in $A := \mathbb{F}_q[T]$. Let $\mathbb{Z}_p$ denote the $p$-adic integers and let $\Omega$ be the completion of an algebraic closure of $K_\infty$. The analogy with characteristic zero is given by $A \leftrightarrow \mathbb{Z}$, $A^+ \leftrightarrow \mathbb{Z}^+$, $K \leftrightarrow \mathbb{Q}$, $K_\infty \leftrightarrow \mathbb{R}$, and $\Omega \leftrightarrow \mathbb{C}$. The Goss zeta function for $\mathbb{F}_q[T]$ is defined as

$$\zeta(z) := \sum_{n \in A^+} n^{-z}$$

where $z$ is taken from $\Omega^* \times \mathbb{Z}_p$. Exponentiation is defined as follows: for a monic polynomial $n$ set $\langle n \rangle := nT^{-\deg n}$; then for $z = (x, y) \in \Omega^* \times \mathbb{Z}_p$ Goss defines

$$n^z := x^{\deg n} \langle n \rangle^y.$$

The term $\langle n \rangle^y$ is well defined since $\langle n \rangle \equiv 1 (\bmod\, T^{-1})$. One draws an analogy between this definition and complex exponentiation when $\mathbb{C}$ is regarded canonically as $\mathbb{R} \times \mathbb{R}$: for positive integers $n$ we have $n^{(x,y)} = (e^x)^{\log n}(e^{i\log n})^y$. Goss showed that by grouping together terms of the same degree $\zeta$ becomes well defined over all $\Omega^* \times \mathbb{Z}_p$:

$$\zeta(z) := \zeta(x,y) := \sum_{m \geq 0} x^{-m} \left( \sum_{n \in A^+, \ \deg n = m} \langle n \rangle^{-y} \right).$$

From the definition of exponentiation we have $n^{(T^k, k)} = n^k$ for any integer $k$. Define $\zeta(k) := \zeta(T^k, k)$ for $k \in \mathbb{Z}$. Thus when $k > 0$ we have $\zeta(k) = \sum_{n \in A^+} n^{-k}$. The values of $\zeta$ on the set of "integers" $\{(T^k, k) : k \in \mathbb{Z}\}$ are analogous to special values of the Riemann zeta function:

- Carlitz showed the following analog to a theorem of Euler. For $k > 0$ we have
$$\zeta((q-1)k) = \frac{B_{(q-1)k} \tilde{\pi}^{(q-1)k}}{g_{(q-1)k}}$$
where $\tilde{\pi} \in \Omega$ is an analog to $2\pi i$, the $B_{(q-1)k} \in K$ are analogs to the even Bernoulli numbers, and the $g_{(q-1)k} \in A$ are analogs to factorials [Car1][Car2].

- Goss showed (in general, not just for the $\mathbb{F}_q[T]$ case) that for $k > 0$ we have $\zeta(-k) \in A$ and $\zeta(-k) = 0$ when $k \equiv 0 \pmod{q-1}$ [Gos1]. In analogy, the Riemann zeta function is rational on the negative integers and zero on the negative evens.

Consult [Gos2] for more interesting properties of Goss' zeta functions including their connection with cyclotomic extensions and with Drinfeld modules.

In Section 2 we prove the following analog to the Riemann Hypothesis. It states that for fixed $y$ the zeros of $\zeta(x, -y)$ are simple and all lie on the same "real line". For a complete explanation of the analog see [Gos2].

**Theorem 1.1.** *Fix $y \in \mathbb{Z}_p$. As a function of $x$, the zeros of $\zeta(x, -y)$ are simple and lie in $K_\infty$. In fact they lie in the subfield $\mathbb{F}_p((T^{-1}))$.*

In [Car3] Carlitz investigated, among other things, the vanishing of the power sums
$$S'_k(N) := \sum_{n \in A^+, \ \deg n = k} n^N$$

for positive integers $N$. He stated that $S'_k(N) \neq 0$ if and only if there exists a $(k+1)$-tuple $(r_0, r_1, \ldots, r_k) \in \mathbb{N}^{k+1}$ whose terms sum to $N$ and satisfy the following two conditions:

(i) there is no carryover of $p$-adic digits in the sum $N = \sum r_j$;

(ii) $r_j > 0$ and $(p^s - 1) \mid r_j$ for $0 \leq j \leq k - 1$.

(Note that (ii) puts no condition on $r_k$.) Let $U_{k+1}(N)$ be the collection of all such $(k+1)$-tuples. Then Carlitz claimed that $S'_k(N) \neq 0$ if and only if $U_{k+1}(N) \neq \emptyset$. His proof went as follows. For a monic $n \in A^+$ set $n = a_0 + a_1 T^1 + \cdots + a_{k-1} T^{k-1} + T^k$. Then

$$S'_k(N) = \sum \binom{N}{r_0, \ldots, r_k} \sum a_0^{r_0} \cdots a_{k-1}^{r_{k-1}} T^{r_1 + 2r_2 + \cdots + kr_k} \quad (1.1)$$

where the outer sum is over all $(k+1)$-tuples $r = (r_0, \ldots, r_k)$ such that $\sum r_j = N$ and the inner sum is over all $(a_0, \ldots, a_{k-1}) \in (\mathbb{F}_q)^k$. Note that the sum $\sum_{a \in \mathbb{F}_q} a^h$ equals $-1$ when $h$ is a positive multiple of $q - 1$, and equals $0$ otherwise. Thus (1.1) becomes

$$S'_k(N) = \sum \binom{N}{r_0, \ldots, r_k} (-1)^k T^{r_1 + 2r_2 + \cdots + kr_k} \quad (1.2)$$

where the sum is over all $(k+1)$-tuples $(r_0, r_1, \ldots, r_k)$ such that $\sum r_j = N$ and also satisfy condition (ii) above. A well known result of Lucas states that the multinomial coefficient

$$\binom{N}{r_0, \ldots, r_k} = \frac{N!}{r_0! \cdots r_k!}$$

is not equivalent to $0 \pmod{p}$ if and only if there is no carryover of $p$-adic digits in the sum $\sum r_j$. So we can take the sum in (1.2) over $U_{k+1}(N)$. Thus if $U_{k+1}(N) = \emptyset$ then $S'_k(N) = 0$. To obtain the converse Carlitz asserted without proof that

> the degree $r_1 + 2r_2 + \cdots + kr_k$ of a monomial appearing
> in (1.2) attains its unique maximum when $(r_k, r_{k-1}, \ldots, r_0)$ (1.3)
> is lexicographically largest among all elements of $U_{k+1}(N)$.

We now present an equivalent but slightly different version of (1.3 ). The changes make our later notation and definitions a bit easier to digest. For an $m$-tuple $X = (X_1, \cdots, X_m) \in \mathbb{N}^m$ we define the *weight* of $X$ as

$$wt(X) := X_1 + 2X_2 + \cdots + mX_m. \quad (1.4)$$

Given a finite subset $W \subset \mathbb{N}^m$, a tuple $O \in W$ is said to be *optimal* in $W$ if $wt(O) \geq wt(X)$ for all $X \in W$. The *greedy* element of $W$ is that tuple $(G_1, \ldots, G_m) \in W$ for which $(G_m, G_{m-1}, \ldots, G_1)$ is largest lexicographically.

A *composition* of $N \in \mathbb{Z}^+$ is a tuple $X = (X_1, \cdots, X_m)$ of *positive* integers that sum to $N$. We say that the $m$-tuple $X$ is a *valid* composition of $N$ if, in addition, there is no carryover of $p$-adic digits in the sum $N = \sum X_j$, and $(p^s - 1) \mid X_j$ for $1 \leq j \leq m - 1$.

Define $V_m(N)$ to be the set of all valid compositions of $N$ of length $m$. Note that
$$V_m(N) = \{(X_1, \cdots, X_m) \in U_m(N) : X_m > 0\} \tag{1.5}$$
Most of this paper is devoted to proving the following theorem.

**Theorem 1.2.** *If $V_m(N)$ is not empty then it contains a unique optimal element. Further, the optimal element is the greedy element of $V_m(N)$.*

Theorem 1.2 implies its counterpart for $U_m(N)$ which we state in Lemma 1.3 below. The proof of Lemma 1.3 uses Proposition 4.6. This proposition states that if $(p^s - 1)$ divides $N$, and $G$ is the greedy element of $V_m(N)$, then

$$(m-1)N < wt(G) \leq mN. \tag{1.6}$$

**Lemma 1.3.** *If $U_m(N)$ is not empty then it contains a unique optimal element. Further, the optimal element is the greedy element of $U_m(N)$.*

**Proof.** First assume that $V_m(N)$ is not empty. Since $V_m(N) \subseteq U_m(N)$ the lemma follows from Theorem 1.2 once we show that $V_m(N)$ contains the greedy and all optimal elements of $U_m(N)$. Note that if $(X_1, \ldots, X_m) \in U_m(N)$, then $X_m \equiv N \pmod{p^s - 1}$ by definition. Thus, by (1.5), if $N$ is not divisible by $(p^s - 1)$ then $U_m(N) = V_m(N)$. So suppose $(p^s - 1)$ divides $N$. Since $V_m(N)$ is not empty its greedy element coincides with that of $U_m(N)$. By (1.5), all tuples in $U_m(N) \backslash V_m(N)$ have the form $(X_1, \ldots, X_{m-1}, 0)$ and have the same weight as the corresponding elements $(X_1, \ldots, X_{m-1})$ in $V_{m-1}(N)$. Now (1.6) implies that the weights of these tuples are strictly less than the weight of the greedy element in $V_m(N)$. Thus if $V_m(N)$ is not empty then it contains all the optimal elements of $U_m(N)$. If on the other hand $V_m(N)$ is empty then $U_m(N)$ is essentially equal to $V_{m-1}(N)$. The result follows. ∎

By comparing (1.2) with the definition of weight (1.4) one see that the degree of a monomial appearing in $S'_k(N)$ has the form $wt(R) - N$ where $R = (r_0, \ldots, r_k) \in U_{k+1}(N)$. Thus Lemma 1.3 is precisely Carlitz's assertion (1.3). In [Car3] Carlitz actually uses (1.3) in two arguments: the one given above and another concerning the vanishing of the sums

$$S_k(N) = \sum_{n \in A, \ \deg n < k} n^N.$$

The next theorem includes a restatement of these results which are now justified by Theorem 1.2. The proof of Part (a) was given above. The proof of Part (b) is similar, see [Car3]. Part (c) follows from Carlitz's assertion (1.3) as was observed in [Tha].

**Theorem 1.4.** *(a) $S'_k(N) \neq 0$ if and only if $U_{k+1}(N) \neq \emptyset$.*
  *(b) $S_k(N) \neq 0$ if and only if $V_k(N) \neq \emptyset$ and $p^s - 1$ divides $N$.*
  *(c) If $S'_k(N) \neq 0$ then its degree is $wt(G) - N$ where $G$ is the greedy element of $U_{k+1}(N)$.*

We now outline Sections 3-7 which make up the proof of Theorem 1.2. In Section 3 we set up our notation, present most of our definitions, and prove some basic results. In Section 4 we give a few more definitions and list several technical results. We end Section 4 with a proof of Theorem 1.2 for the case $s = 1$ (i.e., $q = p$). Sections 5, 6, and 7 constitute the proof for the case $s \geq 2$. Form the set $\Psi$ of all pairs $(m, N)$ such that $V_m(N)$ contains an optimal element which is not the greedy element. Then Theorem 1.2 is equivalent to the statement "$\Psi$ is the empty set." For each positive integer $N$ let $\ell(N)$ be the sum of the $p$-digits of $N$. We assume that $\Psi$ is not empty and choose $(m, N) \in \Psi$ so that $(m, \ell(N))$ is lexicographically minimal. Under this assumption and using our chosen $m$ and $N$, we show that there exists an optimal composition $O$ in $V_m(N)$ whose last part $O_m$ is much smaller than the last part $G_m$ of the greedy composition in $V_m(N)$. The existence of $O$ is established in Section 5. In Section 6 we use $G_m$ and the components of $O$ to carefully construct a third composition $Z \in V_m(N)$. Then we show in Section 7 that $wt(Z)$ is strictly larger than $wt(O)$. This contradicts the fact that $O$ is optimal and so we conclude that $\Psi$ is indeed empty.

## 2. Proof of Theorem 1.1

Let $p$ be prime and set $q = p^s$. In this section we follow [D-V] and use Theorem 1.2 to prove Theorem 1.1: for fixed $y \in \mathbb{Z}_p$ the zeros of $\zeta(x, -y)$, as a function of $x$, are simple and lie in $\mathbb{F}_p((T^{-1}))$.

**Proof.** (of Theorem 1.1) For $(x, y) \in \Omega^* \times \mathbb{Z}_p$ we have

$$\zeta(x, -y) := \sum_{m \geq 0} x^{-m} \left( \sum_{n \in A^+,\ \deg n = m} (T^{-m} n)^y \right).$$

Fix $y \in \mathbb{Z}_p$ and view $\zeta(x, -y)$ as a power series in the variable $x^{-1}$. If $y \in \mathbb{N} \setminus \{0\}$, then $T^{-my} S'_m(y)$ is the coefficient of $x^{-m}$. Further, (1.2) implies these coefficients are in $\mathbb{F}_p((T^{-1}))$. Since $\mathbb{N}$ is dense in $\mathbb{Z}_p$, the coefficients of $\zeta(x, -y)$ are in $\mathbb{F}_p((T^{-1}))$ for any $y$.

Define $v_m(y)$ to be the valuation of the coefficient of $x^{-m}$ in $\zeta(x, -y)$. The *Newton polygon* for $\zeta(x, -y)$ is the lower convex hull in $\mathbb{R}^2$ of the points

$$\{(m, v_m(y)) : m \geq 0\}.$$

Its sides describe the valuations of the zeros of $\zeta(x, -y)$: if the Newton polygon for $\zeta(x, -y)$ has a side of slope $\lambda$ whose projection onto the horizontal axis has length $l$, then $\zeta(x, -y)$ has precisely $l$ zeros with valuation $\lambda$. Here, zeros are counted with multiplicity. When $v_{m-1}(y)$ and $v_m(y)$ are both finite define

$$\lambda_y(m) := v_m(y) - v_{m-1}(y).$$

Then when defined, $\lambda_y(m)$ is the slope of the line segment

$$\text{from } (m-1, v_{m-1}(y)) \text{ to } (m, v_m(y)). \tag{2.1}$$

In Case I below we show that if $\zeta(x, -y)$ is a polynomial of degree $d$ then $\lambda_y$ is both defined and strictly increasing on the set $\{1, 2, \ldots, d\}$; in Case II we show that if $\zeta(x, -y)$ is not a polynomial then $\lambda_y$ is defined and strictly increasing on all of $\mathbb{Z}^+$. In each case we will have shown that these line segments (2.1) are the sides of the Newton polygon for $\zeta(x, -y)$. Since the line segments have horizontal length one, we will have shown that each root of $\zeta(x, -y)$ is simple and lies in $\mathbb{F}_p((T^{-1}))$.

Case I: Assume $y$ is a positive integer. Then $T^{-my}S'_m(y)$ is the coefficient of $x^{-m}$ in $\zeta(x,-y)$. Since $U_{m+1}(y)$ is empty for large enough $m$, Theorem 1.4(a) implies that $\zeta(x,-y)$ is a polynomial of $x^{-1}$. Let $d$ be the degree of $\zeta(x,-y)$. Then $U_k(y)$ is empty for $k > d+1$. Note that if $1 \leq k \leq d+1$ then $U_k(y) \neq \emptyset$: clearly $U_{d+1}(y) \neq \emptyset$; if $(X_1, \ldots, X_{d+1}) \in U_{d+1}(N)$ then for any $1 \leq k \leq d$ we have $(X_1, \ldots, X_k + \cdots + X_{d+1}) \in U_k(N)$. By Theorem 1.4(a) we have

$$v_m(y) = \begin{cases} my - \deg S'_m(y) & \text{if } 0 \leq m \leq d \\ \infty & \text{if } m > d \end{cases}. \tag{2.2}$$

Thus $\lambda_y$ is defined on $\{1, \ldots, d\}$. Further, for $2 \leq m \leq d$ we have

$$\begin{aligned}
\lambda_y(m) - \lambda_y(m-1) &= (v_m(y) - v_{m-1}(y)) - (v_{m-1}(y) - v_{m-2}(y)) \tag{2.3} \\
&= [\deg S'_{m-1}(y) - \deg S'_m(y)] - [\deg S'_{m-2}(y) - \deg S'_{m-1}(y)] \\
&= [\deg S'_{m-1}(y) - \deg S'_{m-2}(y)] - [\deg S'_m(y) - \deg S'_{m-1}(y)].
\end{aligned}$$

Let $F = (F_1, \ldots, F_{m-1})$, $G = (G_1, \ldots, G_m)$, and $H = (H_1, \ldots, H_{m+1})$ be the greedy elements from $U_{m-1}(y)$, $U_m(y)$ and $U_{m+1}(y)$. By Theorem 1.4(c) we have

$$\lambda(m) - \lambda(m-1) = [wt(G) - wt(F)] - [wt(H) - wt(G)]. \tag{2.4}$$

Combine the last two components of $H$ to form $Z := (H_1, \ldots, H_{m-1}, H_m + H_{m+1})$. Then $Z$ is an element of $U_m(y)$ and $wt(Z) = wt(H) - H_{m+1}$. By Lemma 1.3 we have $wt(Z) \leq wt(G)$. Thus from (2.4) we get

$$\begin{aligned}
\lambda(m) - \lambda(m-1) &\geq [wt(Z) - wt(F)] - [wt(H) - wt(Z)] \\
&= [wt(Z) - wt(F)] - H_{m+1}. \tag{2.5}
\end{aligned}$$

Drop the last component of $H$ to form $Y := (H_1, \ldots, H_{m-1}, H_m)$. Since $H$ is the greedy element of $U_{m+1}(y)$, the composition $Y$ is the greedy element of $U_m(y - H_{m+1})$ (Proposition 4.1(a)). Also, we have

$$wt(Y) = wt(Z) - mH_{m+1}. \tag{2.6}$$

Note that $Y$ is also the greedy element in $V_m(y - H_{m+1})$. Since $(q-1) \mid y - H_{m+1}$, Proposition 4.6 implies

$$0 < wt(Y) - (m-1)(y - H_{m+1}). \tag{2.7}$$

Add $H_{m+1}$ to each side of (2.7) and then combine with (2.6) to get

$$\begin{aligned} H_{m+1} &< wt(Y) + mH_{m+1} - (m-1)y \\ &= wt(Z) - (m-1)y. \end{aligned} \tag{2.8}$$

Now from (2.5), first apply the fact that the weight of $F = (F_1, \ldots, F_{m-1})$ is less than or equal to $(m-1)y$ and then apply (2.8) to get

$$\begin{aligned} \lambda_y(m) - \lambda_y(m-1) &\geq [wt(Z) - wt(F)] - H_{m+1} \\ &\geq wt(Z) - (m-1)y - H_{m+1} \\ &> 0. \end{aligned}$$

This completes Case I.

Case II: Assume that $y \in \mathbb{Z}_p - \mathbb{N}$ so that $y = y_0 + y_1 p + \cdots + y_i p^i + \cdots$ has an infinite number of nonzero digits. We need the following easy lemma.

**Lemma 2.1.** *Let $\tilde{y}(t)$ be the sum of the first $t$ terms in the $p$-adic expansion of $y$:*

$$\tilde{y}(t) := \sum_{i=0}^{t} y_i p^i.$$

*Then for any fixed $m$ there exists a positive integer $t'$ such that $U_m(\tilde{y}(t))$ is nonempty for all $t \geq t'$.*

**Proof.** Fix $m$. We construct $t'$ and an element $X = (X_1, \ldots, X_m)$ of $U_m(\tilde{y}(t))$ for any $t \geq t'$. For $0 \leq h \leq s-1$ set

$$z_h := \sum_{i=0}^{\infty} y_{is+h} p^{is+h}.$$

Then $y = z_0 + \cdots + z_{s-1}$. There exists an $h'$ such that $z_{h'}$ has an infinite number of nonzero $p$-adic digits. Form the nondecreasing infinite sequence $\sigma(z_{h'})$ consisting of $y_{is+h'}$ copies of $p^{is+h'}$ for each $i \geq 0$. Set $X_1$ equal to the sum of the first $p^s - 1$ terms of $\sigma(z_{h'})$. Set $X_2$ equal to the next $p^s - 1$ terms and continue in this way up through $X_{m-1}$. Define $t'$ so that $p^{t'}$ is the smallest term in $\sigma(z_{h'})$ which was not assigned. For any fixed $t \geq t'$ set $X_m = \tilde{y}(t) - (X_1 + \cdots + X_{m-1})$. It is now easy to see that $X \in U_m(\tilde{y}(t))$. ∎

Let $\tilde{y}(t)$ be as in the lemma. It is sufficient to show that for any $m$ there exists a $t_m$ such that if $t \geq t_m$ then $v_m(y) = v_m(\tilde{y}(t))$. For in that case $\lambda_y$ is defined for all $m \geq 1$. Further, (2.3) implies

$$\lambda_y(m) - \lambda_y(m-1) = \lambda_{\tilde{y}(t)}(m) - \lambda_{\tilde{y}(t)}(m-1)$$

for $m \geq 2$ and any $t \geq \max\{t_{m-2}, t_{m-1}, t_m\}$. Thus Case II is reduced to Case I.

We show that for all $m \geq 0$ there exists a $t_m$ such that if $t \geq t_m$ then $v_m(y) = v_m(\tilde{y}(t))$. Since $v_0(y) = 0$ for all $y$, we may assume $m \geq 1$. By Lemma 2.1, there exists a positive integer $t'$ such that $U_{m+1}(\tilde{y}(t))$ is nonempty for all $t \geq t'$. For $t \geq t'$ denote the greedy element of $U_{m+1}(\tilde{y}(t))$ by $G^t = (G_1^t, G_2^t, \ldots, G_{m+1}^t)$. Note that for any $X \in U_{m+1}(\tilde{y}(t))$ we have

$$\begin{aligned} wt(X) &= X_1 + \cdots + mX_m + (m+1)X_{m+1} \text{ and} \\ 0 &= (m+1)\tilde{y}(t) - (m+1)X_1 - \cdots - (m+1)X_{m+1}. \end{aligned}$$

By adding we deduce

$$wt(X) - \tilde{y}(t) = m\tilde{y}(t) - (mX_1 + \cdots + 2X_{m-1} + X_m).$$

Thus for $t \geq t'$ we have by Theorem 1.4(c) and (2.2),

$$v_m(\tilde{y}(t)) = m\tilde{y}(t) - \deg S'_m(\tilde{y}(t)) = mG_1^t + \cdots + 2G_{m-1}^t + G_m^t. \tag{2.9}$$

Thus $v_m(\tilde{y}(t))$ depends only on the first $m$ terms of $G^t$. We claim that there exists $t_m$ such that if $t \geq t_m$ then $G^t$ and $G^{t_m}$ differ only in their $(m+1)^{st}$ term. That is, we have $G_i^t = G_i^{t_m}$ for $1 \leq i \leq m$, so that $v_m(\tilde{y}(t)) = v_m(\tilde{y}(t_m))$. To see this, note first that $G^t = (G_1^t, G_2^t, \ldots, G_{m+1}^t) \in U_{m+1}(\tilde{y}(t))$ implies

$$(G_1^t, \ldots, G_m^t, G_{m+1}^t + y_{t+1}p^{t+1}) \in U_{m+1}(\tilde{y}(t+1)). \tag{2.10}$$

Secondly, note that by greediness, $\tilde{y}(t+1) - G_{m+1}^{t+1}$ is minimal in the set

$$\{\tilde{y}(t+1) - X_{m+1} : X \in U_{m+1}(\tilde{y}(t+1))\}. \tag{2.11}$$

Thus (2.10) and (2.11) imply

$$\tilde{y}(t+1) - G_{m+1}^{t+1} \leq \tilde{y}(t+1) - (G_{m+1}^t + y_{t+1}p^{t+1}) = \tilde{y}(t) - G_{m+1}^t.$$

In other words, $\tilde{y}(t) - G^t_{m+1}$ monotonically decreases as $t$ increases. Since $\tilde{y}(t) - G^t_{m+1}$ is bounded below by zero, we can choose $t_m$ large enough so that

$$\tilde{y}(t) - G^t_{m+1} = \tilde{y}(t_m) - G^{t_m}_{m+1} \text{ for all } t \geq t_m.$$

Finally, note that the greedy element of $U_m\left(\tilde{y}(t) - G^t_{m+1}\right)$ is obtained by dropping the last term from the greedy element $G^t = (G^t_1, G^t_2, \ldots, G^t_{m+1})$ of $U_{m+1}(\tilde{y}(t))$ (Proposition 4.1(a)). Thus we have if $t \geq t_m$ then $G^t_i = G^{t_m}_i$ for $1 \leq i \leq m$. Now (2.9) implies we have $v_m(\tilde{y}(t)) = v_m(\tilde{y}(t_m))$ for all $t \geq t_m$. Since $y = \lim_{t\to\infty} \tilde{y}(t)$ we have $v_m(y) = v_m(\lim_{t\to\infty} \tilde{y}(t)) = \lim_{t\to\infty} v_m(\tilde{y}(t)) = v_m(\tilde{y}(t_m))$. ∎

## 3. Notation and definitions

In this section we set up the notation, present some definitions, and state some basic results. We begin with a discussion of the conditions a composition $X = (X_1, \ldots, X_m)$ of $N$ must meet in order to be valid:

(i) there is no carryover of $p$-adic digits in the sum $N = \sum X_j$;
(ii) $X_j$ is a multiple of $p^s - 1$ for $1 \leq j \leq m - 1$.

To deal with condition (i) we treat a given positive integer $N$ as the shortest nondecreasing sequence $\sigma(N)$ of powers of $p$ whose terms sum to $N$. For $N > 0$ define $\deg_p(N)$ to be the exponent of the largest power of $p$ appearing in $\sigma(N)$.

**Example.** Let $p = 3$ and $N = 131$. In base 3 we have $N = 11212_3$. Thus $\sigma(N) = (1, 1, 3, 3^2, 3^2, 3^3, 3^4)$ and $\deg_p(N) = 4$.

We view a sequence as an ordered *multiset*: a set in which an element may appear more than once. A *partition* of a multiset $M$ is a collection of non-empty multisets whose disjoint union is $M$. Thus, a composition $(X_1, \cdots, X_m)$ of $N$ satisfies (i) if and only if $\{\sigma(X_1), \cdots, \sigma(X_m)\}$ is a partition of $\sigma(N)$.

To deal with condition (ii) we use the fact that a number is divisible by $p^s - 1$ if and only if the sum of its $p^s$-adic digits is divisible by $p^s - 1$. In order to keep track of $p$-adic digits when summing $p^s$-adic digits we define the map $\Gamma : \mathbb{N} \to \mathbb{N}^s$ as follows. Given $N \in \mathbb{N}$ with $p$-adic expansion $N = \sum_{j\geq 0} n_j p^j$, define $\Gamma(N)$ to be the column vector $[u_0, u_1, \ldots, u_{s-1}]^t$ where $u_i$ is the sum of all $n_j$ such that $j \equiv i \pmod{s}$ (and the superscript 't' indicates transpose). If $\Gamma(N) = [u_0, u_1, \ldots, u_{s-1}]^t$ then the sum of the $p^s$-adic digits of $N$ is $u_0 + u_1 p + \cdots + u_{s-1} p^{s-1}$.

**Example.** *For $p = 3$ and $s = 2$, we have $\Gamma(11212_3) = [5, 2]^t$.*

Set $\bar{\psi}_0 := [1, p, \ldots, p^{s-1}]^t$ and let $\langle *, * \rangle$ be the standard inner product on $\mathbb{R}^s$. Then a composition $(X_1, \cdots, X_m)$ of $N$ satisfies condition (ii) if and only if

$$(p^s - 1) \mid \langle \bar{\psi}_0, \Gamma(X_j) \rangle \text{ for } 1 \leq j \leq m - 1.$$

Notice that $\Gamma$ is not in general an additive function. However, when $(X_1, \ldots, X_m)$ satisfies condition (i) we have for any $\{j_1, \ldots, j_k\} \subseteq \{1, 2, \ldots, m\}$,

$$\Gamma(X_{j_1}) + \cdots + \Gamma(X_{j_k}) = \Gamma(X_{j_1} + \cdots + X_{j_k}).$$

Define the following partial order: for two vectors $\bar{x} = [x_0, \ldots, x_{s-1}]^t$ and $\bar{y} = [y_0, \ldots, y_{s-1}]^t$ in $\mathbb{R}^s$ we write $\bar{x} \leq \bar{y}$ if and only if $x_i \leq y_i$ for $0 \leq i \leq s - 1$. It is important to note that $\bar{x} < \bar{y}$ means $x_i \leq y_i$ for $0 \leq i \leq s - 1$ with $x_i < y_i$ for at least one $i$.

**Example.** *Set $p = 3$ and $s = 2$. We construct all valid compositions of $11212_3$ with $m = 2$ components. In order to satisfy (i) we partition $\sigma(11212_3)$ into two parts $(\Sigma_1, \Sigma_2)$; then we set $X_i$ equal to the sum of the elements of $\Sigma_i$. In order to satisfy (ii) we take $\Sigma_1$ so that 8 divides $\langle \bar{\psi}_0, \Gamma(X_1) \rangle$. The only vectors $\bar{v} \in \mathbb{N}^2 \setminus \{\bar{0}\}$ such that $\bar{v} < [5, 2]^t$ and 8 divides $\langle \bar{\psi}_0, \bar{v} \rangle = v_0 + 3v_1$ are $[5, 1]^t$ and $[2, 2]^t$. Now it is easy to verify that $V_2(11212_3) = \{(11202_3, 10_3), (11110_3, 102_3), (11011_3, 201_3), (10212_3, 1000_3), (1210_3, 10002_3), (1111_3, 10101), (1012_3, 10200_3)\}$.*

Given a composition $X = (X_1, \ldots, X_m)$ define $\Gamma X$ to be the $s \times m$ matrix with columns $\Gamma(X_1), \ldots, \Gamma(X_m)$.

**Example.** *Set $p = 3$ and $s = 2$. If $X$ is $(11202_3, 10_3)$ or $(10212_3, 1000_3)$ then $\Gamma X = \begin{bmatrix} 5 & 0 \\ 1 & 1 \end{bmatrix}$. For any other $X \in V_2(11212_3)$ above, we have $\Gamma X = \begin{bmatrix} 2 & 3 \\ 2 & 0 \end{bmatrix}$.*

In the proof of Theorem 1.2 we assume there exist $m$ and $N$ such that $V_m(N)$ contains an optimal composition $O$ which is different from the greedy composition. Our goal is to arrive at a contradiction by constructing a composition $Z \in V_m(N)$ such that $wt(Z) > wt(O)$. To construct $Z$ we first construct an $s \times m$ matrix B and then define $Z$ to be optimal among those valid compositions $X$ for which $\Gamma X = B$. We now discuss the construction of $Z$ from B.

For $0 \leq i \leq s - 1$ define $\tau_i(N)$ to be the subsequence of $\sigma(N)$ consisting of all $p^k \in \sigma(N)$ such that $k \equiv i \pmod{s}$.

**Example**. Set $p = 3$, $s = 2$, and $N = 11212_3$. Then $\sigma(N) = (1, 1, 3, 3^2, 3^2, 3^3, 3^4)$, $\Gamma(N) = [5, 2]^t$ and we have $\tau_0(N) = (1, 1, 3^2, 3^2, 3^4)$ and $\tau_1(N) = (3, 3^3)$.

Note that if $\Gamma(N) = [u_0, u_1, \ldots, u_{s-1}]^t$ then $u_i$ is the length of $\tau_i(N)$.

Consider two components $X_i$ and $X_j$, with $i < j$, of a valid composition $X = (X_1, \ldots, X_m)$ of $N$. Suppose $p^k$ is a term in $\tau_h(X_i)$ and $p^l$ is a term of $\tau_h(X_j)$ (so that $k \equiv l \equiv h \pmod{s}$). It is easily seen that

$$X' := (X_1, \ldots, X_i - p^k + p^l, \ldots, X_j - p^l + p^k, \ldots, X_m)$$

is also a valid composition of $N$. We say a valid composition $X = (X_1, \ldots, X_m)$ is $\tau$-*monotonic* if and only if for all $1 \leq i < j \leq m$ and $0 \leq h \leq s - 1$, the largest term of $\tau_h(X_i)$ is no larger than the smallest term of $\tau_h(X_j)$. Equivalently, $X = (X_1, \ldots, X_m)$ is $\tau$-monotonic if and only if for all $0 \leq h \leq s - 1$ the sequence $\tau_h(N)$ is simply the concatenation of the subsequences $\tau_h(X_1), \ldots, \tau_h(X_m)$.

**Lemma 3.1.** *The greedy and all optimal elements of $V_m(N)$ are $\tau$-monotonic.*

**Proof.** If $k > l$ in the definition of $X'$ above, then $X'$ has a larger weight than $X$ and $(X_m, \ldots, X_j - p^l + p^k, \ldots, X_i - p^k + p^l, \ldots, X_1)$ is lexicographically larger than $(X_m, \ldots, X_1)$. ∎

Given an $s \times m$ matrix B, denote by $V_m^B(N)$ the set of all valid compositions $X$ of $N$ such that $\Gamma X = B$.

**Lemma 3.2.** *If $V_m^B(N)$ is not empty then it contains a unique $\tau$-monotonic composition.*

**Proof.** A composition $X = (X_1, \ldots, X_m) \in V_m^B(N)$ is uniquely determined by the sequences $\tau_h(X_j)$. If $X$ is $\tau$-monotonic then for each fixed $h$ the sequences $\tau_h(X_j)$, $1 \leq j \leq m$, are uniquely determined by the sequence $\tau_h(N)$ and the $h^{th}$ row of the matrix B. ∎

**Example**. Set $p = 3$ and $s = 2$. We construct the $\tau$-monotonic composition $(Z_1, Z_2)$ in $V_2^B(11212_3)$ where $B = \begin{bmatrix} 2 & 3 \\ 2 & 0 \end{bmatrix}$. The sequences

$$\tau_0(Z_1) := (1, 1), \quad \tau_0(Z_2) := (3^2, 3^2, 3^4)$$
$$\tau_1(Z_1) := (3, 3^3), \quad \tau_1(Z_2) := \emptyset.$$

are forced by $\tau$-monotonicity. Thus we have $Z_1 := (1+1) + (3+3^3) = 1012_3$ and $Z_2 := 3^2 + 3^2 + 3^4 = 10200_3$.

We now specify the conditions B must meet in order for $V_m^B(N)$ to be non-empty. To this end we use a matrix to characterize the set of column vectors

$$\mathfrak{J} := \{\Gamma(k) : k \text{ is a positive multiple of } p^s - 1\}.$$

Let $\bar{e}_0, \ldots, \bar{e}_{s-1}$ denote the standard basis of column vectors for $\mathbb{R}^s$. Define $\bar{\varepsilon}_i := p\bar{e}_{i-1} - \bar{e}_i$ for $0 \leq i \leq s-1$. Here and for now on indices which should range from 0 to $s-1$ are evaluated modulo $s$: i.e., $a_i = a_{\tilde{i}}$ when $i \equiv \tilde{i} \pmod{s}$. Thus $\bar{\varepsilon}_0 = p\bar{e}_{s-1} - \bar{e}_0$. Define the $s \times s$ matrix $E := [\bar{\varepsilon}_0, \bar{\varepsilon}_1, \ldots, \bar{\varepsilon}_{s-1}]$. For example if $p = 5$ and $s = 3$ then

$$E = \begin{bmatrix} -1 & 5 & 0 \\ 0 & -1 & 5 \\ 5 & 0 & -1 \end{bmatrix}.$$

Below we show that $\mathfrak{J} = (E\mathbb{Z}^s) \cap (\mathbb{N}^s \setminus \{\bar{0}\})$. We will often use the fact that for two arbitrary column vectors $\bar{u} = [u_0, \ldots, u_{s-1}]^t$ and $\bar{a} = [a_0, \ldots, a_{s-1}]^t$, if $\bar{u} = E\bar{a}$ then

$$u_j = pa_{j+1} - a_j.$$

To construct the inverse of $E$, define $R := [\bar{e}_1, \bar{e}_2, \ldots, \bar{e}_{s-1}, \bar{e}_0]$ to be the permutation matrix which rotates the coordinates of a vector to the right: $R\bar{e}_i = \bar{e}_{i+1}$. For $1 \leq i \leq s-1$ define $\bar{\psi}_i := R^i \bar{\psi}_0$. Then

$$\langle \bar{\psi}_i, \bar{\varepsilon}_j \rangle = \begin{cases} p^s - 1 & \text{if } i = j \\ 0 & \text{otherwise} \end{cases}. \tag{3.1}$$

Thus $E^{-1} = (p^s - 1)^{-1} [\bar{\psi}_0, \bar{\psi}_1, \ldots, \bar{\psi}_{s-1}]^t$. We list some simple facts about the vectors $\bar{\psi}_i$.

**Lemma 3.3.** (a) For $k \in \mathbb{N}$, $k \equiv \langle \bar{\psi}_0, \Gamma(k) \rangle \pmod{p^s - 1}$.
  (b) For $\bar{u} \in \mathbb{Z}^s$, $\langle \bar{\psi}_0, R\bar{u} \rangle \equiv p \langle \bar{\psi}_0, \bar{u} \rangle \pmod{p^s - 1}$.
  (c) For $\bar{u} \in \mathbb{Z}^s$, $\langle \bar{\psi}_0, \bar{u} \rangle \equiv p^i \langle \bar{\psi}_i, \bar{u} \rangle \pmod{p^s - 1}$.
  (d) For $\bar{x} \in \mathbb{R}^s$, $\langle \bar{\psi}_i, E\bar{x} \rangle = x_i(p^s - 1)$.

**Proof.** Part (a) is clear. For (b) note that

$$\langle \bar{\psi}_0, R\bar{u} \rangle = p \left( \sum_{i=0}^{s-1} u_i p \right) - u_{s-1}(p^s - 1) = p \langle \bar{\psi}_0, \bar{u} \rangle - u_{s-1}(p^s - 1).$$

Part (c) follows from (b). Part (d) follows from (3.1). ∎

The next lemma gives us a useful characterization of the set $\mathfrak{J}$.

**Lemma 3.4.** *Let $k$ be a positive integer. Then $k$ is a multiple of $p^s - 1$ if and only if $\Gamma(k) = \mathrm{E}\bar{\mathrm{a}}$ for some $\bar{\mathrm{a}} \in \mathbb{Z}^s$. In other words $\mathfrak{J} = (\mathrm{E}\mathbb{Z}^s) \cap (\mathbb{N}^s \setminus \{\bar{0}\})$.*

**Proof.** Suppose $\Gamma(k) = \mathrm{E}\bar{\mathrm{a}}$ for some $\bar{\mathrm{a}} \in \mathbb{Z}^s$. Then by Lemma 3.3(d), $\langle \bar{\psi}_0, \Gamma(k) \rangle = a_0(p^s - 1)$. By Lemma 3.3(a), $k$ is a multiple of $p^s - 1$.

Suppose $k$ is a positive multiple of $p^s - 1$. Set $\bar{\mathrm{x}} := \mathrm{E}^{-1}\Gamma(k)$. It is sufficient to show that $\bar{\mathrm{x}} \in \mathbb{Z}^s$. From the expression for $\mathrm{E}^{-1}$ given above,

$$\bar{\mathrm{x}} = (p^s - 1)^{-1} \left[ \bar{\psi}_0, \bar{\psi}_1, \ldots, \bar{\psi}_{s-1} \right]^{\mathrm{t}} \Gamma(k) = \sum_{i=0}^{s-1} (p^s - 1)^{-1} \langle \bar{\psi}_i, \Gamma(k) \rangle \bar{\mathrm{e}}_i.$$

By Lemma 3.3(a), $p^s - 1$ divides $\langle \bar{\psi}_0, \Gamma(k) \rangle$. Lemma 3.3(c) implies that each $\langle \bar{\psi}_i, \bar{\mathrm{u}} \rangle$ is divisible by $p^s - 1$. Thus $\bar{\mathrm{x}} \in \mathbb{Z}^s$. ∎

Now we state precisely when $V_m^{\mathrm{B}}(N)$ is non-empty.

**Lemma 3.5.** *Let $\mathrm{B} = \left[ \bar{\mathrm{b}}_1, \ldots, \bar{\mathrm{b}}_m \right]$ be an integer matrix. The set $V_m^{\mathrm{B}}(N)$ is nonempty if and only if (a) the columns of $\mathrm{B}$ sum to $\Gamma(N)$ and (b) $\bar{\mathrm{b}}_1, \ldots, \bar{\mathrm{b}}_{m-1}$ are members of $\mathfrak{J}$ and $\bar{\mathrm{b}}_m > \bar{0}$.*

**Proof.** If $X \in V_m^{\mathrm{B}}(N)$ then $\mathrm{B} = \Gamma X$ certainly satisfies (a) and (b). Conversely, suppose that $\mathrm{B} = [b_{i,j}]$ satisfies (a) and (b). We construct the $\tau$-monotonic element of $V_m^{\mathrm{B}}(N)$. By condition (a) one can form subsequences $\Theta_{i,1}, \Theta_{i,2}, \ldots, \Theta_{i,m}$ of $\tau_i(N)$ such that $\Theta_{i,j}$ has length $b_{i,j}$ and such that $\tau_i(N)$ is the concatenation of the sequences $\Theta_{i,1}, \Theta_{i,2}, \ldots, \Theta_{i,m}$ for $0 \leq i \leq s - 1$. For $1 \leq j \leq m$ define $X_j$ to be the sum of the elements of $\bigsqcup_{i=0}^{s-1} \Theta_{i,j}$. Then $\Gamma X = \mathrm{B}$. By construction $\{\sigma(X_1), \cdots, \sigma(X_m)\}$ is a partition of $\sigma(N)$. Condition (b) implies that $X_j$ is a positive multiple of $p^s - 1$ for $1 \leq j \leq m - 1$ and $X_m > 0$. Thus $X \in V_m^{\mathrm{B}}(N)$. ∎

# 4. Some preliminary results and a proof of the $q = p$ case.

This section lists several technical results and ends with a proof of Theorem 1.2 for the case $q = p$. Our first result consists of several observations on how to obtain new optimal (or greedy) compositions from old ones.

**Proposition 4.1.** *Suppose $X = (X_1, \ldots, X_m)$ is an optimal (or the greedy) composition in $V_m(N)$. Then*

*(a) $(X_1, \ldots, X_{m-1})$ is an optimal (or the greedy) composition in $V_{m-1}(N - X_m)$.*

*(b) If $p^k \in \sigma(X_m)$ and $X_m > p^k$ then $(X_1, \ldots, X_{m-1}, X_m - p^k)$ is an optimal (or the greedy) composition in $V_m(N - p^k)$.*

*(c) For any integer $n \geq 0$, the composition $(p^n X_1, \ldots, p^n X_m)$ is an optimal (or the greedy) composition in $V_m(p^n N)$.*

*Furthermore the six statements remain true when the sets $V_*(*)$ are replaced with the sets $U_*(*)$ throughout.*

**Proof.** The twelve statements are all easily proved with similar arguments. The four statements of Part (a) are easiest. We prove Parts (b) and (c) for the sets $V_*(*)$.

For part (b), note that the assumptions $p^k \in \sigma(X_m)$ and $X_m > p^k$ imply that $X^{(b)} := (X_1, \ldots, X_m - p^k)$ is a valid composition in $V_m(N - p^k)$. Suppose $X^{(b)}$ is not optimal (or not greedy). Then there exists a $Y = (Y_1, \ldots, Y_m) \in V_m(N - p^k)$ such that $wt(Y) > wt(X^{(b)})$ (or such that $(Y_m, \ldots, Y_1)$ is lexicographically larger than $(X_m - p^k, \ldots, X_1)$). Clearly $(Y_1, \ldots, Y_m + p^k) \in V_m(N)$. Note that the weight of $(Y_1, \ldots, Y_m + p^k)$ is strictly greater than $wt(X)$ (or that $(Y_m + p^k, \ldots, Y_1)$ is lexicographically larger than $(X_m, \ldots, X_1)$). This is a contradiction.

For part (c), first note that for all $i$ the $i^{th}$ term of $\sigma(p^n N)$ is just $p^n$ times the $i^{th}$ term of $\sigma(N)$. Thus $\{\sigma(p^n X_1), \ldots, \sigma(p^n X_m)\}$ is a partition of $\sigma(p^n N)$. Since $p^s - 1 \mid p^n X_j$ for $j \leq m - 1$ we have

$$X^{(c)} := (p^n X_1, \ldots, p^n X_m) \in V_m(p^n N).$$

Suppose $X^{(c)}$ is not optimal (or not greedy). Then there exists a $Z = (Z_1, \ldots, Z_m) \in V_m(p^n N)$ such that $wt(Z) > wt(X^{(c)})$ (or such that $(Z_m, \ldots, Z_1)$ is lexicographically larger than $(p^n X_m, \ldots, p^n X_1)$). Since each term in $\sigma(p^n N)$ is divisible by $p^n$, we have $\{\sigma(p^{-n} Z_1), \ldots, \sigma(p^{-n} Z_m)\}$ is a partition of $\sigma(N)$. Clearly $p^{-n} Z_j$ is divisible by $p^s - 1$ for $1 \leq j \leq m - 1$. Thus

$$(p^{-n} Z_1, \ldots, p^{-n} Z_m) \in V_m(N).$$

This is a contradiction since the weight of $(p^{-n}Z_1, \ldots, p^{-n}Z_m)$ is strictly greater than $wt(X)$ (or $(p^{-n}Z_m, \ldots, p^{-n}Z_1)$ is lexicographically larger than $(X_m, \ldots, X_1)$). ∎

The remaining results of this section depend on certain subsets of $\mathbb{N}^s$ whose definitions we now motivate. Let $X = (X_1, \ldots, X_m) \in V_m^B(N)$ where $B = [\bar{b}_1, \bar{b}_2, \ldots, \bar{b}_m]$. Suppose that for some $j \leq m-1$ we have $\bar{b}_j = \bar{u}_1 + \bar{u}_2$ with $\bar{u}_1, \bar{u}_2 \in \mathfrak{J}$. Then Lemma 3.5 (applied to the matrix $[\bar{u}_1, \bar{u}_2]$) implies there exists a valid composition $(Y_1, Y_2)$ of $X_j$. Clearly

$$X' := (X_1, \ldots, X_{j-1}, Y_1, X_{j+1}, \ldots, X_m + Y_2)$$

is a valid composition of $N$. Further, it is clear that $wt(X') > wt(X)$ and also that $X'$ is, in our reverse lexicographic order, larger than $X$. Thus $X$ can be neither optimal nor greedy. In other words, if $X = (X_1, \ldots, X_m)$ is either optimal or greedy then for $1 \leq j \leq m-1$ we have $V_2(X_j) = \emptyset$. To take advantage of arguments such as this we define following the subsets of $\mathbb{N}^s$:

$$\begin{aligned} I_m &:= \{\Gamma(k) : k \in \mathbb{N} \text{ and } V_m(k) \neq \emptyset\}; \\ J_m &:= \mathfrak{J} \cap (I_m \setminus I_{m+1}). \end{aligned}$$

We argued above that if $(X_1, \ldots, X_m)$ is either optimal or greedy then $\Gamma(X_j) \in J_1$ for $1 \leq j \leq m-1$.

We have two alternate ways of expressing $I_m$:

$$\begin{aligned} I_m &= \{\bar{u} \in \mathbb{N}^s : \exists \bar{v}_1, \ldots, \bar{v}_{m-1} \in \mathfrak{J} \text{ such that } \bar{u} > \bar{v}_1 + \cdots + \bar{v}_{m-1}\} \quad (4.1) \\ &= \{\bar{u} \in \mathbb{N}^s : \exists \bar{\eta}_1, \ldots, \bar{\eta}_{m-1} \in \mathfrak{J} \text{ such that } \bar{u} > \bar{\eta}_{m-1} > \cdots > \bar{\eta}_1\}. \quad (4.2) \end{aligned}$$

Expression (4.1) follows from Lemma 3.5 by taking

$$B = [\bar{v}_1, \cdots, \bar{v}_{m-1}, \bar{u} - (\bar{v}_1 + \cdots + \bar{v}_{m-1})].$$

Expression (4.2) is gotten from (4.1) by taking $\bar{\eta}_j = \bar{v}_1 + \cdots + \bar{v}_j$. It is clear that

$$I_m \subseteq I_k \text{ whenever } k \leq m.$$

From (4.1) we see that $J_m$ *is the set of all $\bar{u} \in \mathfrak{J}$ such that $\bar{u}$ decomposes into a sum of $m$ elements of $\mathfrak{J}$ but not into a sum of $m+1$ elements of $\mathfrak{J}$.*

Proposition 4.3 below is the seminal result concerning the subsets $I_m$ and $J_m$. Before we prove it we must list some basic facts about the matrices E and R. Recall that $E = [\bar{\varepsilon}_0, \bar{\varepsilon}_1, \ldots, \bar{\varepsilon}_{s-1}]$ where $\bar{\varepsilon}_i = p\,\bar{e}_{i-1} - \bar{e}_i$ and that $R = [\bar{e}_1, \bar{e}_2, \ldots, \bar{e}_{s-1}, \bar{e}_0]$ is the matrix which rotates the coordinates of a vector to the right.

**Lemma 4.2.** *(a)* $\{E^{-1}\bar{x} : \bar{x} > \bar{0}\} \subseteq (\mathbb{R}^+)^s$.
  *(b) If $\bar{x} > \bar{y}$ then $E^{-1}\bar{x} > E^{-1}\bar{y}$.*
  *(c) $RE = ER$.*
  *(d) $R\mathfrak{J} = \mathfrak{J}$.*
  *(e) For any $k \in \mathbb{N}$ we have $R^k \Gamma(N) = \Gamma(p^k N)$.*

**Proof.** Part (a) follows from the fact that all the terms in $E^{-1}$ are positive. Part (b) is equivalent to Part (a). For Part (c) note that $R\bar{\varepsilon}_i = \bar{\varepsilon}_{i+1}$ for all $i$. Thus we have
$$RE = [\bar{\varepsilon}_1, \bar{\varepsilon}_2, \ldots, \bar{\varepsilon}_{s-1}, \bar{\varepsilon}_0] = [E\bar{e}_1, E\bar{e}_2, \ldots, E\bar{e}_{s-1}, E\bar{e}_0] = ER.$$
Now (d) follows from Lemma 3.4 and the fact that R and E commute. For Part (e) let $N = n_0 + n_1 p + \cdots + n_l p^l$ be the $p$-adic expansion of $N$. Set $[u_0, \ldots, u_{s-1}]^t := \Gamma(N)$ and $[v_0, \ldots, v_{s-1}]^t := \Gamma(p^k N)$. Since $p^k N = n_0 p^k + n_1 p^{1+k} + \cdots + n_l p^{l+k}$ is the $p$-adic expansion of $p^k N$ we have $v_{i+k} = u_i$. ∎

**Proposition 4.3.** *For $m \geq 1$ we have,*
  *(a) $I_m = \{E\bar{x} \in \mathbb{N}^s : m - 1 < \min\{x_0, \ldots, x_{s-1}\}\}$, and*
  *(b) $J_m = \{E\bar{a} \in \mathbb{N}^s : m = \min\{a_0, \ldots, a_{s-1}\}\}$*
  *where $\bar{x} = [x_0, \ldots, x_{s-1}]^t$ and $\bar{a} = [a_0, \ldots, a_{s-1}]^t$ are taken from $\mathbb{R}^s$.*

**Proof.** First we show that (a) implies (b). Assume (a) and let $\bar{a} \in \mathbb{R}^s$. Suppose first that $E\bar{a} \in J_m$. By definition $E\bar{a} \in \mathfrak{J} \cap (I_m \backslash I_{m+1})$. Since $E\bar{a} \in (I_m \backslash I_{m+1})$, Part (a) implies
$$m - 1 < \min\{a_0, \ldots, a_{s-1}\} \leq m.$$
Since $E\bar{a} \in \mathfrak{J}$ we have $\bar{a} \in \mathbb{Z}^s$. Thus $m = \min\{a_0, a_1, \ldots, a_{s-1}\}$ and therefore
$$J_m \subseteq \{E\bar{a} \in \mathbb{N}^s : m = \min\{a_0, , \ldots, a_{s-1}\}\}.$$
For the other inclusion suppose $E\bar{a} \in \mathbb{N}^s$ and $m = \min\{a_0, \ldots, a_{s-1}\}$. Then (a) implies $E\bar{a} \in (I_m \backslash I_{m+1})$. To show $E\bar{a} \in \mathfrak{J}$ we show by reverse induction that $\bar{a} \in \mathbb{Z}^s$. Suppose $a_h = m$ and suppose further that $a_{i+1}$ is an integer for some $i + 1 \leq h$. Set $[w_0, \ldots, w_{s-1}]^t := E\bar{a}$. Then
$$w_i = pa_{i+1} - a_i.$$
Since $w_i$ and $a_{i+1}$ are integers, so is $a_i$. Thus $\bar{a} \in \mathbb{Z}^s$. Therefore $E\bar{a} \in \mathfrak{J} \cap (I_m \backslash I_{m+1}) = J_m$ and (a) implies (b).

We prove Part (a). Suppose $\bar{u} \in I_m$. Set $\bar{x} := E^{-1}\bar{u}$. By (4.1) there exist vectors $\bar{v}_1, \ldots, \bar{v}_m$ in $\mathfrak{J}$ such that $\bar{u} > \bar{v}_1 + \cdots + \bar{v}_{m-1}$. For $1 \leq j \leq m-1$, set $\bar{b}_j := E^{-1}\bar{v}_j$. Lemma 4.2(b) implies $\bar{x} > (\bar{b}_1 + \cdots + \bar{b}_{m-1})$. By Lemma 4.2(a) and the definition of $\mathfrak{J}$, each $\bar{b}_j := [b_{0,j}, \ldots, b_{s-1,j}]^t$ is in $(\mathbb{N}\backslash\{0\})^s$. Thus, for $0 \leq i \leq s-1$ we have $x_i > b_{i,1} + \cdots + b_{i,m-1} \geq m-1$.

Now suppose that $\bar{u} = E\bar{x} \in \mathbb{N}^s$ with $m - 1 < \min\{x_0, \ldots, x_{s-1}\}$. Let $k$ be such that $x_k = \min\{x_0, \ldots, x_{s-1}\}$. Let $\lfloor x \rfloor$ denote the largest integer $\leq x$. Then $\lfloor x_k \rfloor \geq m - 1$. By (4.2) our proof will be complete once we show that there exist vectors $\bar{\eta}_1, \ldots, \bar{\eta}_{\lfloor x_k \rfloor} \in \mathfrak{J}$ such that $\bar{u} > \bar{\eta}_{\lfloor x_k \rfloor} > \cdots > \bar{\eta}_1$. We lose no generality by assuming $k = 0$: if $k \neq 0$, we construct $\bar{\eta}_1, \ldots, \bar{\eta}_{\lfloor x_k \rfloor}$ using $R^{s-k}\bar{u}$ in place of $\bar{u}$. Then $\bar{u} > R^k\bar{\eta}_{\lfloor x_k \rfloor} > \cdots > R^k\bar{\eta}_1$ and the $R^k\bar{\eta}_i$ are in $\mathfrak{J}$ by Lemma 4.2(d). Note that $\bar{u}$ satisfies the conditions on $\bar{v}$ stated in Lemma 4.4 below. The proof is completed by recursively applying Lemma 4.4 to get $\bar{\eta}_{\lfloor x_k \rfloor}, \bar{\eta}_{\lfloor x_k \rfloor - 1}, \ldots, \bar{\eta}_1$.

**Lemma 4.4.** *Suppose that $\bar{v} = E\bar{c} \in \mathbb{N}^s$ with $1 < c_0 = \min\{c_0, \ldots, c_{s-1}\}$. Then there exists a vector $\bar{w} = E\bar{d} \in \mathfrak{J}$ with $\bar{w} < \bar{v}$ and such that $c_0 - 1 \leq d_0 = \min\{d_0, \ldots, d_{s-1}\}$.*

**Proof.** We define the terms of $\bar{d} = [d_0, \ldots, d_{s-1}]^t$ inductively as follows where $\lceil x \rceil$ is the smallest integer $\geq x$.

$$\text{Set } d_0 := \lceil c_0 \rceil - 1.$$
$$\text{Set } d_i := \min\{\lceil c_i \rceil - 1, pd_{i+1}\} \text{ for } i = s-1, s-2, \ldots, 1.$$

We claim that
$$d_0 = \min\{d_0, \ldots, d_{s-1}\}. \tag{4.3}$$

Since $c_0 \leq c_{s-1}$ we have $d_{s-1} = \min\{\lceil c_{s-1}\rceil - 1, pd_0\} \geq d_0$. Proceeding by reverse induction, suppose $d_{i+1} \geq d_0$ for some $i+1 \leq s-1$. Then we have $d_i = \min\{\lceil c_i \rceil - 1, pd_{i+1}\} \geq d_0$. Thus (4.3) is justified and therefore $c_0 - 1 \leq d_0 = \min\{d_0, \ldots, d_{s-1}\}$.

Note that now we have
$$d_i = \min\{\lceil c_i \rceil - 1, pd_{i+1}\} \text{ for } all \; i. \tag{4.4}$$

With $\bar{w} := [w_0, \ldots, w_{s-1}]^t := E\bar{d}$ we have
$$\begin{aligned} w_i &= pd_{i+1} - d_i \\ &= \max\{pd_{i+1} - \lceil c_i \rceil + 1, 0\}. \end{aligned}$$

This implies $\bar{w} \in \mathfrak{J}$.

We show that $\bar{w} < \bar{v}$. Fix $i$ and assume $w_i > 0$. Then $w_i = pd_{i+1} - \lceil c_i \rceil + 1$. Since $v_i = pc_{i+1} - c_i$ we have

$$v_i - w_i = p(c_{i+1} - d_{i+1}) + \lceil c_i \rceil - c_i - 1.$$

By (4.4) we have $c_{i+1} > d_{i+1}$. Since $\lceil c_i \rceil \geq c_i$ we have $v_i - w_i > -1$. Since $v_i - w_i$ is an integer we must have $w_i \leq v_i$. Further, $\bar{w} \neq \bar{v}$ since $d_0 < c_0$. Thus $\bar{w} < \bar{v}$ as desired. ∎ ∎

Define
$$\bar{\rho} := [p-1, \ldots, p-1]^{\mathrm{t}}.$$

Our next result is a simple lemma used by the proof of Proposition 4.6 below.

**Lemma 4.5.** *Let $\bar{u}$ be an element of $J_m$ with $m \geq 1$. Suppose $\bar{v} \in \mathbb{N}^s$ with $\bar{v} \leq \bar{\rho}$ and $\bar{v} < \bar{u}$. Then $\bar{u} - \bar{v}$ is an element of $I_m \cup J_{m-1}$ (where $J_0 := \emptyset$).*

**Proof.** Since $I_1 = \mathbb{N}^s \setminus \{\bar{0}\}$, we may assume $m > 1$. Set $\bar{a} := [a_0, \ldots, a_{s-1}]^{\mathrm{t}} := \mathrm{E}^{-1}\bar{u}$ and $\bar{x} := [x_0, \ldots, x_{s-1}]^{\mathrm{t}} := \mathrm{E}^{-1}\bar{v}$. Note that $\mathrm{E}^{-1}\bar{\rho} = [1, \ldots, 1]^{\mathrm{t}}$. Lemma 4.2(b) implies $x_i \leq 1$ for all $i$. Thus by Proposition 4.3(b), we have $m - 1 \leq \min\{a_0 - x_0, \ldots, a_{s-1} - x_{s-1}\}$. By Proposition 4.3 again, we have $\bar{u} - \bar{v} \in I_m \cup J_{m-1}$. ∎

**Proposition 4.6.** *Fix $N$ and $m$ such that $p^s - 1$ divides $N$ and $V_m(N)$ is not empty. If $X$ is either an optimal or the greedy composition in $V_m(N)$ then*

$$(m-1)N < wt(X) \leq mN.$$

**Proof.** Note that $wt(X) \leq mN$ is true for any $X$ in $V_m(N)$. Our proof of the other inequality is by induction on $m$. The result is trivially true when $m = 1$. Suppose it is true for $m < m_0$. Fix $m = m_0 > 1$. Let

$$N = n_0 + n_1 p + \cdots + n_k p^k$$

be the $p$-adic expansion of $N$. Since $N$ is divisible by $p^s - 1$ we have $\bar{u} := \Gamma(N) \in J_h$ for some $h \geq m$. Set $\bar{v} := \Gamma(n_k p^k)$. As $1 \leq n_k \leq p-1$ we have $\bar{v} \leq \bar{\rho}$. Also $\bar{v} < \bar{u}$ since $m > 1$ rules out the possibility that $N = n_k p^k$. By Lemma 4.5, $\bar{u} - \bar{v}$ is an element of $I_h \cup J_{h-1}$. Since $h \geq m$ and $\bar{u} - \bar{v} = \Gamma(N - n_k p^k)$ we have $\Gamma(N - n_k p^k) \in I_m \cup J_{m-1}$.

Suppose that $\Gamma(N - n_k p^k) \in I_m$. Then $V_m(N - n_k p^k)$ is not empty. Let $(Y_1, \ldots, Y_m)$ be any composition in $V_{m-1}(N - n_k p^k)$ and set

$$Y' := (Y_1, \ldots, Y_m + n_k p^k).$$

If $\Gamma(N - n_k p^k) \in J_{m-1}$ choose $Y$ from $V_{m-1}(N - n_k p^k)$ and set

$$Y' := (Y_1, \ldots, Y_{m-1}, n_k p^k).$$

In either case $Y' \in V_m(N)$ and $Y'_m \geq n_k p^k$.

Define $G := (G_1, \ldots, G_m)$ to be the greedy composition in $V_m(N)$. Greediness implies $G_m \geq Y'_m$. Thus

$$G_m \geq n_k p^k > N/2. \tag{4.5}$$

Set $Z := (G_1, \ldots, G_{m-1})$. Then $Z$ is the greedy composition in $V_{m-1}(N - G_m)$ (Proposition 4.1(a)). Further, we have

$$wt(G) = wt(Z) + mG_m. \tag{4.6}$$

Our induction hypothesis implies

$$(m-2)(N - G_m) < wt(Z). \tag{4.7}$$

Combining (4.6) and (4.7) we get

$$wt(G) > (m-2)(N - G_m) + mG_m = (m-2)N + 2G_m.$$

Now (4.5) implies $(m-1)N < wt(G)$. If $O$ is an optimal composition in $V_m(N)$ then $wt(O)$ must satisfy the same inequality. ∎

**Proposition 4.7.** *If $X = (X_1, \ldots, X_m)$ is either an optimal or the greedy composition in $V_m(N)$ then $\Gamma(N - X_m) \in J_{m-1}$.*

**Proof.** We show the contrapositive: if $\Gamma(N - X_m) \notin J_{m-1}$ then $X$ is neither greedy nor optimal.

Since $p^s - 1$ divides $N - X_m$ we have $\Gamma(N - X_m) \in \mathfrak{J} \cap I_{m-1}$. Assume $\Gamma(N - X_m) \notin J_{m-1}$. Then by the definitions we have $\Gamma(N - X_m) \in I_m$ and thus $V_m(N - X_m) \neq \emptyset$. Let $Y := (Y_1, \ldots, Y_m)$ be any composition in $V_m(N - X_m)$. Then $Y' := (Y_1, \ldots, Y_m + X_m)$ is a valid composition in $V_m(N)$. Thus $X$ cannot be greedy since $(Y_m + X_m, \ldots, Y_1)$ is lexicographically larger than $(X_m, \ldots, X_1)$.

Set $X' := (X_1, \ldots, X_{m-1}) \in V_{m-1}(N - X_m)$. Then $wt(X) = wt(X') + mX_m$. Let $G := (G_1, \ldots, G_m)$ be the greedy composition in $V_m(N - X_m)$. Set

$$Z := (G_1, \ldots, G_{m-1}, G_m + X_m).$$

Then $Z \in V_m(N)$ and $wt(Z) = wt(G) + mX_m$. Since $p^s - 1$ divides $N - X_m$, Proposition 4.6 implies

$$wt(X') \leq (m-1)(N - X_m) < wt(G).$$

Adding $mX_m$ across these inequalities gives $wt(X) < wt(Z)$. Thus $X$ cannot be optimal. ∎

Our next result concerns the following subsets of $J_m$:

$$J_m^i := \{\mathrm{E}\bar{a} \in \mathbb{N}^s : a_i = m = \min\{a_0, \ldots, a_{s-1}\}\}.$$

**Proposition 4.8.** *If $X = (X_1, \ldots, X_m)$ is either an optimal or the greedy composition in $V_m(N)$ then there exists an $i$ such that $\Gamma(X_j) \in J_1^i$ for $1 \leq j \leq m-1$.*

**Proof.** Note that

$$\begin{aligned} \Gamma(X_1) + \cdots + \Gamma(X_{m-1}) &= \Gamma(X_1 + \cdots + X_{m-1}) \\ &= \Gamma(N - X_m). \end{aligned}$$

We argued above that each $\Gamma(X_j)$ must be in $J_1$. Since $\Gamma(N - X_m) \in J_{m-1}$ by Proposition 4.7, the result is implied by Proposition 4.3. ∎

**Proposition 4.9.** *Suppose $X = (X_1, \ldots, X_m)$ is either greedy or optimal in $V_m(N)$ and suppose $p^k$ is any element of $\sigma(X_m)$. Then $V_m(N - p^k)$ is empty if and only if $X_m = p^k$.*

**Proof.** If $X_m > p^k$ then $(X_1, \ldots, X_m - p^k)$ would be an element of $V_m(N - p^k)$. Thus, if $V_m(N - p^k)$ is empty then $X_m = p^k$. Suppose that $X_m = p^k$. Then $\Gamma(N - p^k) \in J_{m-1}$ by Proposition 4.7. Therefore $\Gamma(N - p^k) \notin I_m$ which is another way of saying $V_m(N - p^k)$ is empty. ∎

Now we prove Theorem 1.2 for the case $q = p$. Recall that $\Psi$ is the set of all pairs $(m, N)$ such that $V_m(N)$ contains an optimal element that is not the greedy element.

**Proposition 4.10.** *If $s = 1$ then $\Psi$ is empty.*

**Proof.** Fix $m$ and $N$ so that $V_m(N)$ is not empty. Let $G$ be the greedy and $O$ any optimal element of $V_m(N)$. Note that since $s = 1$ we have $J_1 = \{p-1\}$. Thus $\Gamma G = \Gamma O = \mathrm{B}$ where B is the $1 \times m$ matrix

$$[\underbrace{p-1, \ldots, p-1}_{m-1}, \Gamma(N) - (m-1)(p-1)].$$

(In this case $\Gamma(N) = \ell(N) =$ the sum of the $p$-digits of $N$.) By Lemmas 3.1 and 3.2 we have $G = O =$ the $\tau$-monotonic element of $V_m^{\mathrm{B}}(N)$. ∎

## 5. If $\Psi \neq \emptyset$ then there exists an optimal $O\ldots$

Sections 5, 6, and 7 constitute the proof of Theorem 1.2 for the case $s \geq 2$. In these sections we assume that $\Psi$ is not empty in order to derive a contradiction. The purpose of the present section is to prove the following proposition.

**Proposition 5.1.** *Suppose $s \geq 2$. If $\Psi$ is not empty then there exists an $m > 2$, an $N$, and an optimal element $(O_1, \ldots, O_m) \in V_m(N)$ which satisfies the properties listed below. Here $(G_1, \ldots, G_m)$ is the greedy element of $V_m(N)$.*
 (a) $\Gamma(O_j) \in J_1^0$ *for* $1 \leq j \leq m-1$.
 (b) $\deg_p(O_m) < \deg_p(G_m)$.
 (c) $O_m = p^k$ *for some* $k \in \mathbb{N}$.
 (d) $\langle \bar{\psi}_0, \Gamma(G_m) \rangle = \langle \bar{\psi}_0, \Gamma(O_m) \rangle$ *yet* $\Gamma(G_m) \neq \Gamma(O_m)$.

**Proof.** Choose $(m, K) \in \Psi$ so that $(m, \ell(K))$ is lexicographically minimal. Let $O' = (O'_1, \ldots, O'_m)$ be any optimal element of $V_m(K)$ that *differs* from the greedy element $G' = (G'_1, \ldots, G'_m)$. By Proposition 4.8, there exists an $h$ such that $\Gamma(O'_j) \in J_1^h$ for $1 \leq j \leq m-1$. Define $N := p^{s-h}K$ and define

$$O := (O_1, \ldots, O_m) := (p^{s-h}O'_1, \ldots, p^{s-h}O'_m).$$

Lemma 4.2(e) implies Part (a):

$$\Gamma(O_j) = \Gamma(p^{s-h}O'_j) \in J_1^0 \text{ for } 1 \leq j \leq m-1.$$

By Proposition 4.1(c), $O$ is an optimal and
$$G := (G_1, \ldots, G_m) = (p^{s-h}G'_1, \ldots, p^{s-h}G'_m)$$
is the greedy element of $V_m(N)$. Clearly $O \neq G$. Thus $(m, N) \in \Psi$. Further, since $\ell(N) = \ell(K)$, we have that $(m, \ell(N))$ is lexicographically minimal. Thus we have the following lemma.

**Lemma 5.2.** *If $m' < m$ or if $m' = m$ and $\ell(N') < \ell(N)$ then either $V_{m'}(N')$ is empty or its greedy element is its unique optimal element.*

It is obvious that $m > 1$. We argue that $m > 2$. If $G = (G_1, G_2)$ is the greedy element of $V_2(N)$ then $G_2 > X_2$ where $X = (X_1, X_2)$ is any other element of $V_2(N)$. But then
$$wt(G) - N = G_2 > X_2 = wt(X) - N.$$
Thus the greedy element of $V_2(N)$ is the unique optimal element. This implies $m > 2$.

Set
$$\begin{aligned}
\bar{o} &= [o_0, \ldots, o_{s-1}]^t = \Gamma(O_m), \\
\bar{g} &= [g_0, \ldots, g_{s-1}]^t = \Gamma(G_m), \text{ and} \\
\bar{u} &= [u_0, \ldots, u_{s-1}]^t = \Gamma(N).
\end{aligned}$$

Define integers $\beta$ and $\gamma$ by
$$\begin{aligned}
\beta &\equiv \deg_p(O_m) (\text{mod } s), \text{ and} \\
\gamma &\equiv \deg_p(G_m) (\text{mod } s)
\end{aligned}$$
where $0 \leq \beta, \gamma \leq s - 1$.

We show $\deg_p(G_m) > \deg_p(O_m)$. By greediness, we have $\deg_p(G_m) \geq \deg_p(O_m)$. Suppose that $\deg_p(G_m) = \deg_p(O_m)$. Then $g_\beta$ and $o_\beta$ would both be positive. But this contradicts following lemma.

**Lemma 5.3.** *For each $h$ with $0 \leq h \leq s - 1$, either $g_h = 0$ or $o_h = 0$.*

**Proof.** Suppose there exists an $h$ such that $g_h > 0$ and $o_h > 0$. By Lemma 3.1 both $O$ and $G$ are $\tau$-monotonic. Thus $\sigma(G_m)$ and $\sigma(O_m)$ both contain the largest element, say $p^n$, of $\tau_h(N)$.

Suppose first that $V_m(N - p^n)$ is empty. Then $O_m = G_m = p^n$ by Proposition 4.9. By Proposition 4.1(a), $(O_1, \ldots, O_{m-1})$ is optimal and $(G_1, \ldots, G_{m-1})$ is greedy in $V_{m-1}(N - p^n)$. But then Lemma 5.2 implies $(O_1, \ldots, O_{m-1}) = (G_1, \ldots, G_{m-1})$ and so $O = G$. This contradicts our choice of $O$.

Now suppose $V_m(N - p^n)$ is not empty. Then $O_m$ and $G_m$ are both $> p^n$ by Proposition 4.9. Proposition 4.1(b) implies $(O_1, \ldots, O_m - p^n)$ is optimal and $(G_1, \ldots, G_m - p^n)$ is greedy in $V_m(N - p^n)$. Again Lemma 5.2 implies the contradiction that $O = G$. ∎

Thus we have Part (b).

For Parts (c) and (d) we need Lemma 5.5. Lemma 5.5 needs Lemma 5.4 below.

**Lemma 5.4.** *For any $x$ such that $o_x > 0$ we have*

$$\bar{u} - \bar{e}_\gamma - \bar{e}_x \notin J_{m-1} \cup I_m.$$

**Proof.** Since $o_x > 0$, Lemma 5.3 implies $x \neq \gamma$. Suppose on the contrary that $\bar{u} - \bar{e}_\gamma - \bar{e}_x \in J_{m-1} \cup I_m$. Then there exist vectors $\bar{v}_1, \ldots, \bar{v}_{m-1}$ in $\mathfrak{J}$ such that

$$\bar{v}_1 + \ldots + \bar{v}_{m-1} \leq \bar{u} - \bar{e}_\gamma - \bar{e}_x. \tag{5.1}$$

Let $p^c$ be the largest element in $\tau_x(N)$. Then we have $\Gamma(N - p^c) = \bar{u} - \bar{e}_x$. Since $\bar{v}_1 + \ldots + \bar{v}_{m-1} < \bar{u} - \bar{e}_x$ we have $\bar{0} < \bar{v}_m := (\bar{u} - \bar{e}_x) - (\bar{v}_1 + \ldots + \bar{v}_{m-1})$. Proposition 3.5, applied to the matrix $B := [\bar{v}_1, \ldots, \bar{v}_m]$, implies that $V_m^B(N - p^c)$ is not empty. Let $X$ be the $\tau$-monotonic element of $V_m^B(N - p^c)$. Note that (5.1) implies that the $\gamma$-th coordinate of $\bar{v}_m$ is positive. Thus we have $\deg_p(X_m) = \deg_p(G_m)$.

Since $o_x > 0$ and $O$ is $\tau$-monotonic (Lemma 3.1) we have $p^c \in \sigma(O_m)$. Thus $O_m > p^c$. Proposition 4.1(b) implies that $Y := (O_1, \ldots, O_{m-1}, O_m - p^c)$ is optimal in $V_m(N - p^c)$. But then Lemma 5.2 implies that $Y$ is also greedy in $V_m(N - p^c)$. This is a contradiction since

$$\deg_p(Y_m) \leq \deg_p(O_m) < \deg_p(G_m) = \deg_p(X_m). \blacksquare$$

The Parts (c) and (d) will follow from the next lemma. Recall that $\bar{\psi}_0 = [1, p, \ldots, p^{s-1}]^t$ and $\bar{\psi}_i = R^i \bar{\psi}_0$ for $0 \leq i \leq s - 1$ where the matrix $R$ rotates a vector's coordinates to the right.

**Lemma 5.5.** *Let $x$ be such that $o_x > 0$ and set*
$$\bar{v} := [v_0, \ldots, v_{s-1}]^t := E^{-1}(\bar{u} - \bar{e}_x - \bar{e}_\gamma).$$
*Then there exists an $h$ such that $v_h < m - 1$. Further, for any such $h$ we have*
(i) $\langle \bar{\psi}_h, \bar{g} \rangle = \langle \bar{\psi}_h, \bar{o} \rangle$,
(ii) $\langle \bar{\psi}_h, \bar{o} \rangle < \langle \bar{\psi}_h, \bar{e}_x \rangle + \langle \bar{\psi}_h, \bar{e}_\gamma \rangle$, *and*
(iii) $\langle \bar{\psi}_h, \bar{e}_\gamma \rangle < \langle \bar{\psi}_h, \bar{e}_x \rangle$.

**Proof.** By Lemma 5.4 we have $\bar{u} - \bar{e}_x - \bar{e}_\gamma \notin J_{m-1} \cup I_m$. Thus by Proposition 4.3 we have $\min\{v_0, \ldots, v_{s-1}\} < m - 1$. Fix $h$ such that $v_h < m - 1$.

First we prove Part (i). Set $\bar{a} := E^{-1}(\bar{u} - \bar{o})$ and $\bar{b} := E^{-1}(\bar{u} - \bar{g})$. By Proposition 4.7 both $\bar{u} - \bar{o}$ and $\bar{u} - \bar{g}$ are elements of $J_{m-1}$. Proposition 4.3(b) implies $a_h \geq m - 1$ and $b_h \geq m - 1$. We show that
$$a_h = m - 1 = b_h.$$

First we show $a_h = m - 1$. Note that since $o_x > 0$ we have $\bar{u} - \bar{o} - \bar{e}_\gamma \leq \bar{u} - \bar{e}_x - \bar{e}_\gamma$. Thus
$$\langle \bar{\psi}_h, \bar{u} - \bar{o} - \bar{e}_\gamma \rangle \leq \langle \bar{\psi}_h, \bar{u} - \bar{e}_x - \bar{e}_\gamma \rangle. \tag{5.2}$$
From Lemma 3.3(d) we have $\langle \bar{\psi}_h, \bar{u} - \bar{o} \rangle = a_h(p^s - 1)$ and $\langle \bar{\psi}_h, \bar{u} - \bar{e}_x - \bar{e}_\gamma \rangle = v_h(p^s - 1)$. Evaluating (5.2) gives
$$a_h(p^s - 1) - p^{\gamma'} \leq v_h(p^s - 1)$$
where $p^{\gamma'} = \langle \bar{\psi}_h, \bar{e}_\gamma \rangle$. Thus we deduce
$$a_h - v_h \leq \frac{p^{\gamma'}}{p^s - 1}.$$

Since $v_h < m - 1 \leq a_h$ and $p^{\gamma'} < p^s - 1$, we have $0 < a_h - v_h < 1$. Since $a_h$ is an integer we must have $a_h = m - 1$ (and $m - 2 < v_h < m - 1$). The argument that $b_h = m - 1$ is completely analogous. Since $g_\gamma > 0$ we have $\bar{u} - \bar{g} - \bar{e}_x \leq \bar{u} - \bar{e}_x - \bar{e}_\gamma$. Thus
$$\langle \bar{\psi}_h, \bar{u} - \bar{g} - \bar{e}_x \rangle \leq \langle \bar{\psi}_h, \bar{u} - \bar{e}_x - \bar{e}_\gamma \rangle.$$
Set $p^{x'} = \langle \bar{\psi}_{i_x}, \bar{e}_x \rangle$. As above we get
$$b_h - v_h \leq \frac{p^{x'}}{p^s - 1}.$$

Thus, just as $a_h = m - 1$, we have that $b_h = m - 1$. Lemma 3.3(d) implies

$$\langle \bar{\psi}_h, \bar{u} - \bar{g} \rangle = (m-1)(p^s - 1) = \langle \bar{\psi}_h, \bar{u} - \bar{o} \rangle. \tag{5.3}$$

Thus $\langle \bar{\psi}_h, \bar{g} \rangle = \langle \bar{\psi}_h, \bar{o} \rangle$ as desired.

Now for Part (ii). By (5.3) and Lemma 3.3(d) we have

$$\begin{aligned} \langle \bar{\psi}_h, \bar{u} - \bar{o} \rangle &= (m-1)(p^s - 1) \\ &> v_h(p^s - 1) \\ &= \langle \bar{\psi}_h, \bar{u} - \bar{e}_x - \bar{e}_\gamma \rangle. \end{aligned}$$

From this we deduce

$$\langle \bar{\psi}_h, \bar{o} \rangle < \langle \bar{\psi}_h, \bar{e}_x \rangle + \langle \bar{\psi}_h, \bar{e}_\gamma \rangle. \tag{5.4}$$

Part (iii) is more difficult. Note that since $\gamma \neq x$ we have $\langle \bar{\psi}_h, \bar{e}_\gamma \rangle \neq \langle \bar{\psi}_h, \bar{e}_x \rangle$. To establish $\langle \bar{\psi}_h, \bar{e}_\gamma \rangle < \langle \bar{\psi}_h, \bar{e}_x \rangle$ we show that if $\langle \bar{\psi}_h, \bar{e}_\gamma \rangle > \langle \bar{\psi}_h, \bar{e}_x \rangle$ then $O$ cannot be optimal.

Assume $\langle \bar{\psi}_h, \bar{e}_\gamma \rangle > \langle \bar{\psi}_h, \bar{e}_x \rangle$. Recall that $\gamma \equiv \deg_p(G_m)$. Thus $\bar{e}_\gamma \leq \bar{g}$ and therefore $\langle \bar{\psi}_h, \bar{e}_\gamma \rangle \leq \langle \bar{\psi}_h, \bar{g} \rangle$. Part (i) now implies

$$\langle \bar{\psi}_h, \bar{e}_\gamma \rangle \leq \langle \bar{\psi}_h, \bar{o} \rangle. \tag{5.5}$$

Since $\langle \bar{\psi}_h, \bar{e}_x \rangle < \langle \bar{\psi}_h, \bar{e}_\gamma \rangle$ we must have

$$\langle \bar{\psi}_h, \bar{e}_i \rangle < \langle \bar{\psi}_h, \bar{e}_\gamma \rangle \text{ for all } i \text{ such that } o_i > 0. \tag{5.6}$$

Statement (5.6) follows because if $o_i > 0$ and $i \neq x$ then we have $\langle \bar{\psi}_h, \bar{e}_i \rangle < \langle \bar{\psi}_h, \bar{e}_x \rangle + \langle \bar{\psi}_h, \bar{e}_\gamma \rangle$ by (5.4). But since $\langle \bar{\psi}_h, \bar{e}_i \rangle$, $\langle \bar{\psi}_h, \bar{e}_x \rangle$, and $\langle \bar{\psi}_h, \bar{e}_\gamma \rangle$ are three distinct powers of $p$ and $\langle \bar{\psi}_h, \bar{e}_x \rangle < \langle \bar{\psi}_h, \bar{e}_\gamma \rangle$, we must have $\langle \bar{\psi}_h, \bar{e}_i \rangle < \langle \bar{\psi}_h, \bar{e}_\gamma \rangle$.

Rewrite (5.5) by replacing $\langle \bar{\psi}_h, \bar{o} \rangle$ with its expansion as a sum of powers of $p$:

$$\langle \bar{\psi}_h, \bar{e}_\gamma \rangle \leq \sum_{i,\ o_i > 0} \sum_{j=1}^{o_i} \langle \bar{\psi}_h, \bar{e}_i \rangle. \tag{5.7}$$

By (5.6) each $\langle \bar{\psi}_h, \bar{e}_i \rangle$ appearing on the right hand side of (5.7) divides $\langle \bar{\psi}_h, \bar{e}_\gamma \rangle$. Thus there must be some subset of the $\langle \bar{\psi}_h, \bar{e}_i \rangle$'s which sum to precisely $\langle \bar{\psi}_h, \bar{e}_\gamma \rangle$. In other words, there must exist a vector $\bar{w} \in \mathbb{N}^s$ such that $\bar{w} \leq \bar{o}$ and $\langle \bar{\psi}_h, \bar{w} \rangle =$ $\langle \bar{\psi}_h, \bar{e}_\gamma \rangle$. In yet other words, there must be a positive integer $W$ such that $\sigma(W)$ is a subsequence of $\sigma(O_m)$, and $\Gamma(W) = \bar{w}$, and

$$\langle \bar{\psi}_h, \Gamma(W) \rangle = \langle \bar{\psi}_h, \Gamma(p^{\deg_p(G_m)}) \rangle. \tag{5.8}$$

Using Lemma 3.3 it is easily argued that (5.8) implies $W \equiv p^{\deg_p(G_m)} \pmod{p^s - 1}$. We know that $p^{\deg_p(G_m)}$ is a term in $\sigma(O_j)$ for some $j < m$. Define

$$X := (O_1, \ldots, O_j - p^{\deg_p(G_m)} + W, \ldots, O_m - W + p^{\deg_p(G_m)}).$$

Since congruences $\pmod{p^s - 1}$ are preserved we have $X \in V_m(N)$. Since $\sigma(W)$ is a subsequence of $\sigma(O_m)$ we have $W < p^{\deg_p(G_m)}$ since, as was shown above, $\deg_p(O_m) < \deg_p(G_m)$. Thus $wt(X) > wt(O)$ which is a contradiction. Therefore our supposition that $\langle \bar{\psi}_h, \bar{e}_x \rangle < \langle \bar{\psi}_h, \bar{e}_\gamma \rangle$ must be false. ∎

We use Lemma 5.5 to show Part (c) of the proposition: that $\bar{o}$ is the standard unit vector $\bar{e}_\beta$ where $\beta \equiv \deg_p(O_m) \pmod s$. By Lemma 5.5 (ii) and (iii) there exists $h_\beta$ such that

$$\langle \bar{\psi}_{h_\beta}, \bar{o} \rangle < \langle \bar{\psi}_{h_\beta}, \bar{e}_\beta \rangle + \langle \bar{\psi}_{h_\beta}, \bar{e}_\gamma \rangle \quad \text{and} \tag{5.9}$$

$$\langle \bar{\psi}_{h_\beta}, \bar{e}_\gamma \rangle < \langle \bar{\psi}_{h_\beta}, \bar{e}_\beta \rangle. \tag{5.10}$$

Note that (5.9) and (5.10) imply that $\langle \bar{\psi}_{h_\beta}, \bar{o} \rangle < 2 \langle \bar{\psi}_{h_\beta}, \bar{e}_\beta \rangle$. Thus

$$o_\beta = 1. \tag{5.11}$$

If $s = 2$ then $\bar{o} = \bar{e}_\beta$ now follows from Lemma 5.3 since $\gamma \neq \beta$.

Suppose $s \geq 3$. Before continuing we make a simple observation about the vectors $\bar{\psi}_i$. Fix integers $a, b, c$ such that $0 \leq a < b < c \leq s - 1$. Let $0 \leq i \leq s - 1$ and set $p^{a'} := \langle \bar{\psi}_i, \bar{e}_a \rangle$, $p^{b'} := \langle \bar{\psi}_i, \bar{e}_b \rangle$ and $p^{c'} := \langle \bar{\psi}_i, \bar{e}_c \rangle$. Then it is easily seen that:

$$\begin{aligned} &\text{if } i \leq a \text{ or } c < i && \text{then } p^{a'} < p^{b'} < p^{c'}; \\ &\text{if } a < i \leq b && \text{then } p^{b'} < p^{c'} < p^{a'}; \\ &\text{if } b < i \leq c && \text{then } p^{c'} < p^{a'} < p^{b'}. \end{aligned} \tag{5.12}$$

We continue with our proof that $\bar{o} = \bar{e}_\beta$ when $s \geq 3$. Suppose that in addition to $\beta$ there exists another index $\alpha \neq \beta$ such that $o_\alpha > 0$. Applying Lemma 5.5 again we have that there exists $h_\alpha$ such that

$$\langle \bar{\psi}_{h_\alpha}, \bar{o} \rangle < \langle \bar{\psi}_{h_\alpha}, \bar{e}_\alpha \rangle + \langle \bar{\psi}_{h_\alpha}, \bar{e}_\gamma \rangle \quad \text{and} \tag{5.13}$$

$$\left\langle \bar{\psi}_{h_\alpha}, \bar{e}_\gamma \right\rangle < \left\langle \bar{\psi}_{h_\alpha}, \bar{e}_\alpha \right\rangle.$$

Since $o_\alpha > 0$ and $o_\beta > 0$ with $\alpha \neq \beta$, (5.9) implies that we must have $\left\langle \bar{\psi}_{h_\beta}, \bar{e}_\alpha \right\rangle < \left\langle \bar{\psi}_{h_\beta}, \bar{e}_\gamma \right\rangle$. Similarly, (5.13) implies we must have $\left\langle \bar{\psi}_{h_\alpha}, \bar{e}_\beta \right\rangle < \left\langle \bar{\psi}_{h_\alpha}, \bar{e}_\gamma \right\rangle$. Now it follows that we have

$$\left\langle \bar{\psi}_{h_\beta}, \bar{e}_\alpha \right\rangle < \left\langle \bar{\psi}_{h_\beta}, \bar{e}_\gamma \right\rangle < \left\langle \bar{\psi}_{h_\beta}, \bar{e}_\beta \right\rangle \quad \text{and} \tag{5.14}$$

$$\left\langle \bar{\psi}_{h_\alpha}, \bar{e}_\beta \right\rangle < \left\langle \bar{\psi}_{h_\alpha}, \bar{e}_\gamma \right\rangle < \left\langle \bar{\psi}_{h_\alpha}, \bar{e}_\alpha \right\rangle. \tag{5.15}$$

This is a contradiction since (5.12) implies that we cannot have both (5.14) and (5.15) simultaneously true. Thus no $\alpha$ exists. This and (5.11) implies $\bar{o} = \bar{e}_\beta$. This completes the proof of Part (c) of the proposition.

Now for Part (d). By Lemma 5.3 we have $\Gamma(G_m) \neq \Gamma(O_m)$. To show $\left\langle \bar{\psi}_0, \Gamma(G_m) \right\rangle = \left\langle \bar{\psi}_0, \Gamma(O_m) \right\rangle$ we use Lemma 5.5(i) with $x = \beta$. It is sufficient to show that $v_0 < m - 1$. Note that by Part (c) we have $\bar{u} - \bar{e}_\beta = \Gamma(N) - \Gamma(O_m) = \Gamma(O_1) + \cdots + \Gamma(O_{m-1})$. Thus

$$\bar{v} = E^{-1}(\bar{u} - \bar{e}_\beta - \bar{e}_\gamma) = \left( \sum_{j=1}^{m-1} E^{-1} \Gamma(O_j) \right) - E^{-1} \bar{e}_\gamma.$$

By Part (a), $\Gamma(O_j) \in J_1^0$. Thus $v_0 = (m-1) - p^\gamma/(p^s - 1)$. ∎

## 6. Constructing $Z \in V_m(N)$.

For this section we assume $s \geq 2$ and that $\Psi$ is not empty. We use the notation of Section 5: $O = (O_1, \ldots, O_m) \in V_m(N)$ where $O$, $m$, and $N$ have been chosen so that $O$ has the properties listed in Proposition 5.1; $\beta \equiv \deg_p(O_m) \pmod{s}$, $0 \leq \beta \leq s - 1$; $G = (G_1, \ldots, G_m)$ is the greedy element of $V_m(N)$; $\bar{g} = [g_0, \ldots, g_{s-1}]^t = \Gamma(G_m)$; $\bar{u} = [u_0, \ldots, u_{s-1}]^t = \Gamma(N)$.

Our ultimate goal is to arrive at a contradiction by showing that for some $Z \in V_m(N)$ we have $wt(Z) > wt(O)$. In this section we construct $Z$. We construct a new matrix B from the matrix $\Gamma O$ and the vector $\bar{g}$. Then we let $Z$ be the $\tau$-monotonic element of $V_m^B(N)$. Because the construction of B is complicated we outline it here. As a first step we form the matrix $[\theta_{i,j}] := [\bar{\theta}_1, \ldots, \bar{\theta}_m]$ whose columns are the partial sums of the columns of $\Gamma O$:

$$\bar{\theta}_j := \sum_{i=1}^{j} \Gamma(O_i), \text{ for } 1 \leq j \leq m.$$

Second we take $E^{-1}$ of the result:

$$[w_{i,j}] := [\bar{w}_1, \ldots, \bar{w}_m] := E^{-1}[\theta_{i,j}].$$

Next, using Lemma 6.1 below, we perturb the rows of $[w_{i,j}]$ to obtain the matrix $[d_{i,j}]$. Finally, we reverse the first two steps: set $[\delta_{i,j}] = [\bar{\delta}_1, \ldots, \bar{\delta}_m] = E[d_{i,j}]$; then set $\bar{b}_1 = \bar{\delta}_1$, $\bar{b}_2 = \bar{\delta}_2 - \bar{\delta}_1, \ldots, \bar{b}_m = \bar{\delta}_m - \bar{\delta}_{m-1}$ to obtain $B = [\bar{b}_1, \ldots, \bar{b}_m]$.

To obtain $[d_{i,j}]$ we only perturb those rows of $[w_{i,j}]$ with indices between $\alpha$ and $\beta$ where $\alpha$ is given by Lemma 6.1 below.

**Lemma 6.1.** *There exists an index $\alpha$, $0 \leq \alpha < \beta$ such that*
  *(a) if $g_i > 0$ then $\alpha \leq i < \beta$.*
  *(b) $w_{i,m-1} > m-1$ for $\alpha < i \leq \beta$, and*
  *(c) $w_{\alpha,m-1} = m-1$.*

**Proof.** Set $\bar{v} := \bar{u} - \bar{g}$ and $\bar{c} := E^{-1}\bar{v}$. By Proposition 4.7 we have $\bar{v} \in J_{m-1}$. Thus the terms of $\bar{c}$ are all $\geq m-1$ (Proposition 4.3). Set $\bar{a} := E^{-1}(\bar{g} - \bar{e}_\beta)$. By Proposition 5.1 we have $\Gamma(O_{m,}) = \bar{e}_\beta$ and so $\bar{u} - \bar{e}_\beta = \bar{\theta}_{m-1}$. Thus

$$\bar{c} + \bar{a} = E^{-1}(\bar{u} - \bar{e}_\beta) = E^{-1}(\bar{\theta}_{m-1}) = \bar{w}_{m-1}. \tag{6.1}$$

We show that $\bar{a}$ is in $\mathbb{N}^s$. Note that (6.1) implies $\bar{a} \in \mathbb{Z}^s$ since $\bar{c}$ and $\bar{w}_{m-1}$ are both in $\mathbb{N}^s$. By Proposition 5.1(d) we have $\langle \bar{\psi}_0, \bar{g} - \bar{e}_\beta \rangle = 0$. Thus $a_0 = 0$ (Lemma 3.3(d)). Proceeding by induction, suppose $a_{h-1} \geq 0$. Note that since $\bar{g} := E\bar{a} + \bar{e}_\beta$ we have

$$g_{h-1} = \begin{cases} pa_h - a_{h-1} & \text{if } h-1 \neq \beta \\ pa_h - a_{h-1} + 1 & \text{if } h-1 = \beta \end{cases}. \tag{6.2}$$

Since $g_{h-1} \geq 0$ and $p \geq 2$ we have $a_h \geq 0$. Thus $\bar{a} \in \mathbb{N}^s$.

By Proposition 5.1(d) we have $\langle \bar{\psi}_0, \bar{g} \rangle = \langle \bar{\psi}_0, \bar{e}_\beta \rangle = p^\beta$ yet $\bar{g} \neq \bar{e}_\beta$. Thus we have $g_i = 0$ for $\beta \leq i \leq s-1$. Suppose $\beta < s-1$. Since $a_0 = 0$ and $0 = g_{s-1} = a_0 p - a_{s-1}$ we have $a_{s-1} = 0$. Now an easy induction argument implies that $a_i = 0$ for $\beta < i \leq s-1$. Thus if $\beta < s-1$ then $a_{\beta+1} = 0$. Also, if $\beta = s-1$ then $a_{\beta+1} = a_0 = 0$ as noted above. With $h-1 = \beta$, (6.2) implies $a_\beta = 1$. Let $\alpha'$ be smallest such that $a_{\alpha'+1} > 0$. Then $\alpha' < \beta$. Proceeding by induction, suppose $a_{h-1} > 0$ where $\alpha' + 1 \leq h-1 < \beta$. Since $g_{h-1} \geq 0$, (6.2) implies $a_h > 0$. We conclude that

$$\bar{a} \in \mathbb{N}^s \text{ and } a_i > 0 \text{ if and only if } \alpha' < i \leq \beta. \tag{6.3}$$

Now (6.1) and (6.3) imply

$$w_{i,m-1} > m - 1 \text{ for } \alpha' < i \leq \beta. \tag{6.4}$$

Further, (6.3) and (6.2) imply

$$\text{if } g_i > 0 \text{ then } \alpha' \leq i < \beta. \tag{6.5}$$

Note that Proposition 5.1(a) implies $\bar{\theta}_{m-1} \in J^0_{m-1}$. Thus $w_{0,m-1} = m - 1$. Define $\alpha$ to be the maximum integer $j \leq \alpha'$ such that $w_{j,m-1} = m - 1$. Thus we have Part (c). Parts (a) and (b) follow from (6.5) and (6.4). ∎

We now use $[w_{i,j}]$ to define $[d_{i,j}] = [\bar{d}_1, \ldots, \bar{d}_m]$. First define $\bar{d}_m := E^{-1}\Gamma(N)$. We recursively define the remaining elements of $[d_{i,j}]$. Define

$$d_{i,m-1} := \begin{cases} w_{i,m-1} & \text{if } 0 \leq i \leq \alpha \text{ or } \beta < i \leq s - 1 \\ \min\{w_{i,m-1} - 1, pd_{i+1,m-1}\} & \text{for } i = \beta, \beta - 1, \ldots, \alpha + 1 \end{cases}$$

and for $j = m - 2, m - 3, \ldots, 1$, define

$$d_{i,j} := \begin{cases} w_{i,j} & \text{if } 0 \leq i \leq \alpha \text{ or } \beta < i \leq s - 1 \\ \min\{d_{i,j+1} - 1, pd_{i+1,j}\} & \text{for } i = \beta, \beta - 1, \ldots, \alpha + 1 \end{cases}.$$

Note that for $j \leq m - 1$, since the vectors $\bar{\theta}_j$ are members of $\mathfrak{J}$, we have that the entries of the vectors $\bar{w}_j$ are integers. Thus the definition implies that the $\bar{d}_j$ are integer vectors. Now we define the integer matrix

$$[\delta_{i,j}] := [\bar{\delta}_1, \ldots, \bar{\delta}_m] := E[d_{i,j}]$$

Many of the arguments of Section 7 depend on the following result which will not only allow us to define $Z$ but will also describe its structure in comparison to $O$.

**Lemma 6.2.** *Let* $[u_1, \ldots, u_m]^t := \Gamma(N)$. *Then*
 (a) $\delta_{h,m-1} = u_h$                           *for* $0 \leq h < \alpha$ *or* $\beta \leq h \leq s - 1$,
 (b) $0 \leq \delta_{h,m-1} \leq \max\{u_h - (p-1), 0\}$ *for* $\alpha < h < \beta$,
 (c) $0 \leq \delta_{\alpha,m-1} \leq u_\alpha - p$.
*For* $1 \leq j \leq m - 2$ *we have*
 (d) $\delta_{h,j} = \theta_{h,j}$                            *for* $0 \leq h < \alpha$ *or* $\beta < h \leq s - 1$,
 (e) $0 \leq \delta_{h,j} \leq \max\{\delta_{h,j+1} - (p-1), 0\}$ *for* $\alpha \leq h \leq \beta$.

**Proof.** We start with Part (a). For $0 \leq h < \alpha$ or $\beta < h \leq s - 1$ we have

$$\begin{aligned}
\delta_{h,m-1} &= pd_{h+1,m-1} - d_{h,m-1} \\
&= pw_{h+1,m-1} - w_{h,m-1} \\
&= \theta_{h,m-1}.
\end{aligned}$$

Since $\bar{\theta}_{m-1} = \bar{u} - \bar{e}_\beta$ we have $\delta_{h,m-1} = u_h$ for these indices. Note that we also have

$$\begin{aligned}
\delta_{\beta,m-1} &= pd_{\beta+1,m-1} - d_{\beta,m-1} \\
&= \max\{pw_{\beta+1,m-1} - w_{\beta,m-1} + 1, 0\} \\
&= \max\{\theta_{\beta,m-1} + 1, 0\} \\
&= \theta_{\beta,m-1} + 1 \\
&= u_\beta.
\end{aligned}$$

This completes Part (a).

For Part (b) fix $h$, $\alpha < h < \beta$. Then

$$\begin{aligned}
\delta_{h,m-1} &= pd_{h+1,m-1} - d_{h,m-1} \\
&= \max\{pd_{h+1,m-1} - w_{h,m-1} + 1, 0\}.
\end{aligned}$$

Since $d_{h+1,m-1} \leq w_{h+1,m-1} - 1$ we have

$$\begin{aligned}
0 \leq \delta_{h,m-1} &\leq \max\{pw_{h+1,m-1} - w_{h,m-1} - (p-1), 0\} \\
&= \max\{\theta_{h,m-1} - (p-1), 0\} \\
&= \max\{u_h - (p-1), 0\}.
\end{aligned}$$

So we have Part (b).

For Part (c) we must show $0 \leq \delta_{\alpha,m-1} \leq u_\alpha - p$. We have $\delta_{\alpha,m-1} = pd_{\alpha+1,m-1} - d_{\alpha,m-1}$. Since $d_{\alpha,m-1} = w_{\alpha,m-1}$ and $d_{\alpha+1,m-1} \leq w_{\alpha+1,m-1} - 1$,

$$\begin{aligned}
\delta_{\alpha,m-1} &= pd_{\alpha+1,m-1} - w_{\alpha,m-1} & (6.6) \\
&\leq p(w_{\alpha+1,m-1} - 1) - w_{\alpha,m-1} \\
&= \theta_{\alpha,m-1} - p \\
&= u_i - p.
\end{aligned}$$

We still must show $0 \leq \delta_{\alpha,m-1}$. Since $w_{\alpha,m-1} = m-1$ (Lemma 6.1(c)), we see from (6.6) that it is sufficient to show $d_{\alpha+1,m-1} \geq m-1$. Since $\bar{\theta}_{m-1} \in J_{m-1}$, Proposition 4.3 implies $d_{\beta+1,m-1} = w_{\beta+1,m-1} \geq m-1$. We proceed by reverse induction. Suppose we have $d_{h+1,m-1} \geq m-1$, for some $\alpha < h+1 \leq \beta+1$. Then

$$d_{h,m-1} = \min\{w_{h,m-1} - 1, pd_{h+1,m-1}\}.$$

Since $w_{h,m-1} - 1 \geq m-1$ (Lemma 6.1(b)), we have $d_{h,m-1} \geq m-1$. Thus, by induction

$$d_{h,m-1} \geq m-1 \text{ for } \alpha \leq h \leq \beta, \tag{6.7}$$

and in particular, $d_{\alpha+1,m-1} \geq m-1$.

For Part (d) fix $1 \leq j \leq m-2$. For $0 \leq h < \alpha$ or $\beta < h \leq s-1$ we have

$$\delta_{h,j} = pd_{h+1,j} - d_{h,j} = pw_{h+1,j} - w_{h,j} = \theta_{h,j}.$$

For Part (e) we must show $0 \leq \delta_{h,j} \leq \max\{\delta_{h,j+1} - (p-1), 0\}$ for $\alpha \leq h \leq \beta$ and $1 \leq j \leq m-2$. For these $j$'s for and $\alpha < h \leq \beta$ we have

$$\begin{aligned} \delta_{h,j} &= pd_{h+1,j} - d_{h,j} \\ &= \max\{pd_{h+1,j} - d_{h,j+1} + 1, 0\}. \end{aligned} \tag{6.8}$$

If $h$ is *strictly* less than $\beta$ then $d_{h+1,j} \leq d_{h+1,j+1} - 1$ by definition. Thus

$$0 \leq \delta_{h,j} \leq \max\{\delta_{h,j+1} - (p-1), 0\} \text{ for } \alpha < h < \beta.$$

We deal with the indices $\alpha$ and $\beta$ separately starting with $\beta$. Note that $\mathrm{E}(\bar{\mathrm{w}}_{n+1} - \bar{\mathrm{w}}_n) = \bar{\theta}_{n+1} - \bar{\theta}_n = \Gamma(O_{n+1})$. Since $\Gamma(O_n) \in J_1$ for $1 \leq n \leq m-1$, Proposition 4.3 implies

$$w_{h,n+1} - w_{h,n} \geq 1 \text{ for all } h \text{ and } 1 \leq n \leq m-2. \tag{6.9}$$

By definition $d_{\beta+1,n} = w_{\beta+1,n}$ for each $n$. Thus

$$d_{\beta+1,j} \leq d_{\beta+1,j+1} - 1. \tag{6.10}$$

Combining (6.8) and (6.10) gives

$$0 \leq \delta_{\beta,j} \leq \max\{\delta_{\beta,j+1} - (p-1), 0\}.$$

Now for the index $\alpha$. By Lemma 6.1 we have $w_{\alpha,m-1} = m-1$. Thus (6.9) implies $w_{\alpha,n} = n$ for $1 \le n \le m-1$. By definition $d_{\alpha,n} = w_{\alpha,n}$ for $1 \le n \le m-1$. Thus for $j \le m-2$ we have

$$\begin{aligned} \delta_{\alpha,j} &= pd_{\alpha+1,j} - d_{\alpha,j} \\ &\le p(d_{\alpha+1,j+1} - 1) - (d_{\alpha,j+1} - 1) \\ &= \delta_{\alpha,j+1} - (p-1). \end{aligned} \qquad (6.11)$$

We still must show $\delta_{\alpha,j} \ge 0$. By (6.11) and since $d_{\alpha,j} = j$, it is sufficient to show $d_{h,j} \ge j$ for $\alpha < h \le \beta$. We proceed by reverse induction on $j$ and then on $h$: our basis step is (6.7). Suppose $d_{h,j+1} \ge j+1$ for $\alpha < h \le \beta$ and $j+1 \le m-1$ then

$$d_{\beta,j} = \min\{d_{\beta,j+1} - 1, pd_{\beta+1,j}\} \ge j$$

since $d_{\beta+1,j} = w_{\beta+1,j} \ge j$ (Proposition 4.3). Suppose for some $\alpha < h+1 \le \beta$ we have $d_{h+1,j} \ge j$. Then

$$d_{h,j} = \min\{d_{h,j+1} - 1, pd_{h+1,j}\} \ge j.$$

Thus $d_{h,j} \ge j$ for all $\alpha < h \le \beta$ and $1 \le j \le m-2$. ∎

From the various parts of Lemma 6.2 we have $\bar{0} \le \bar{\delta}_1 \le \ldots \le \bar{\delta}_{m-1} < \bar{\delta}_m$. For $0 \le j \le m-1$, since $\bar{\theta}_j \in J_j^0$ we have $w_{0,j} = j$. Since $d_{0,j} = w_{0,j}$, we have

$$\bar{0} < \bar{\delta}_1 < \ldots < \bar{\delta}_{m-1} < \bar{\delta}_m.$$

For $j \le m-1$, since $\bar{d}_j$ is an integer vector and $\bar{\delta}_j > \bar{0}$, we have $\bar{\delta}_j \in \mathfrak{J}$. Define $\bar{b}_1 := \bar{\delta}_1, \bar{b}_2 := \bar{\delta}_2 - \bar{\delta}_1, \ldots, \bar{b}_m := \bar{\delta}_m - \bar{\delta}_{m-1}$. Then for $j \le m-1$, we have $\bar{b}_j \in \mathfrak{J}$. Now define $B := [\bar{b}_1, \ldots, \bar{b}_m]$. Then the columns of B sum to $\Gamma(N)$. Lemma 3.5 implies that $V_m^B(N)$ is not empty. Finally we can define $Z := (Z_1, \ldots, Z_m)$ to be the $\tau$-monotonic element of $V_m^B(N)$.

## 7. The composition $O$ is not optimal

We continue to assume $s \ge 2$ and that $\Psi$ is not empty. We use the notation of Sections 5 and 6. Thus $m$, $N$, and $O = (O_1, \ldots, O_m) \in V_m(N)$ are chosen so that $O$ has the properties listed in Proposition 5.1. Further, $Z \in V_m(N)$ is as constructed in Section 6. For $1 \leq j \leq m$ set

$$\begin{aligned} \tilde{Z}_j &:= Z_1 + Z_2 + \cdots + Z_j, \text{ and} \\ \tilde{O}_j &:= O_1 + O_2 + \cdots + O_j. \end{aligned}$$

Then we have $\Gamma(\tilde{Z}_j) = \bar{\delta}_j$ and $\Gamma(\tilde{O}_j) = \bar{\theta}_j$ where the matrices $[\delta_{i,j}] = [\bar{\delta}_1, \ldots, \bar{\delta}_m]$ and $[\theta_{i,j}] = [\bar{\theta}_1, \ldots, \bar{\theta}_m]$ are as in Lemma 6.2. In this section we complete the proof of Theorem 1.2 by showing that $wt(Z) > wt(O)$.

It is routine to verify that

$$\begin{aligned} wt(Z) &= mN - (\tilde{Z}_1 + \cdots + \tilde{Z}_{m-1}) \text{ and} \\ wt(O) &= mN - (\tilde{O}_1 + \cdots + \tilde{O}_{m-1}). \end{aligned}$$

Subtracting we get

$$wt(Z) - wt(O) = \sum_{j=1}^{m-1} (\tilde{O}_j - \tilde{Z}_j).$$

Note that $\tilde{Z}_{m-1} = N - Z_m$ and $\tilde{O}_{m-1} = N - O_m$. By Proposition 5.1(c) there exists a positive integer $k$ such that $O_m = p^k$. Thus $\tilde{O}_{m-1} - \tilde{Z}_{m-1} = Z_m - p^k$ and so we have

$$wt(Z) - wt(O) = Z_m - \left( p^k + \sum_{j=1}^{m-2} (\tilde{Z}_j - \tilde{O}_j) \right). \tag{7.1}$$

Since $Z$ is $\tau$-monotonic, we have that $\tau_h(N)$ is the concatenation of the sequences $\tau_h(Z_1), \ldots, \tau_h(Z_m)$. For $1 \leq j \leq m$, since $\tilde{Z}_j = Z_1 + \cdots + Z_j$, we have that $\tau_h(\tilde{Z}_j)$ is the concatenation of the first $j$ sequences $\tau_h(Z_1), \ldots, \tau_h(Z_j)$. Note also that since $\Gamma(\tilde{Z}_j) = \bar{\delta}_j$, the length of $\tau_h(\tilde{Z}_j)$ is $\delta_{h,j}$. For $0 \leq h \leq s - 1$ and $1 \leq i \leq u_h$ define $\tau_{h,i}$ to be the $i^{th}$ term of $\tau_h(N)$. Then, for example,

$$\begin{aligned} \tau_h(\tilde{Z}_j) &= (\tau_{h,1}, \tau_{h,2}, \ldots, \tau_{h,\delta_{h,j}}) \text{ and} \\ \tau_h(Z_m) &= (\tau_{h,\delta_{h,m-1}+1}, \tau_{h,\delta_{h,m-1}+2}, \ldots, \tau_{h,u_h}) \end{aligned}$$

Define $\tau_{h,0} := 0$ for each $h$. This is done so that facts such as

$$\tilde{Z}_j = \sum_{h=0}^{s-1} \sum_{i=0}^{\delta_{h,j}} \tau_{h,i} \quad \text{and} \quad \tilde{O}_j = \sum_{h=0}^{s-1} \sum_{i=0}^{\theta_{h,j}} \tau_{h,i} \tag{7.2}$$

may be unambiguously expressed even when some of the $\delta_{h,j}$ or $\theta_{h,j}$ happen to be zero.

Our next result gives an upper bound on part of expression (7.1).

**Lemma 7.1.** *We have*

$$\sum_{j=1}^{m-2}(\tilde{Z}_j - \tilde{O}_j) \leq p^k + \sum_{h=\alpha}^{\beta-1}\tau_{h,\delta_{h,m-1}}.$$

*where $\alpha$ and $\beta$ are as in Section 6.*

**Proof.** Combining the expressions of (7.2) we have for all $1 \leq j \leq m-2$,

$$\tilde{Z}_j - \tilde{O}_j = \sum_{h=0}^{s-1}\left(\sum_{i=0}^{\delta_{h,j}}\tau_{h,i} - \sum_{i'=0}^{\theta_{h,j}}\tau_{h,i'}\right).$$

By Lemma 6.2(d) we have

$$\tilde{Z}_j - \tilde{O}_j = \sum_{h=\alpha}^{\beta}\left(\sum_{i=0}^{\delta_{h,j}}\tau_{h,i} - \sum_{i'=0}^{\theta_{h,j}}\tau_{h,i'}\right)$$

$$\leq \sum_{h=\alpha}^{\beta}\sum_{i=0}^{\delta_{h,j}}\tau_{h,i}. \tag{7.3}$$

Recall that for each $0 \leq h \leq s-1$ we have $\tau_{h,i} = p^{h+ns}$ for some integer $n \geq 0$. If $\tau_{h,\delta_{h,j}} = p^{h+ts}$, then

$$\sum_{i=0}^{\delta_{h,j}}\tau_{h,i} < p^{h+ts+1}$$

since $\tau_{h,\delta_{h,j}}$ is the largest summand and any particular power of $p$ can appear at most $p-1$ times in the sequence $\tau_h(N)$. By Lemma 6.2(e), if $\delta_{h,j} > 0$ then $\delta_{h,j+1} - \delta_{h,j} \geq p-1$. This implies there are at least $p-1$ terms $x$ of $\tau_h(N)$ such that $\tau_{h,\delta_{h,j}} < x \leq \tau_{h,\delta_{h,j+1}}$. Hence, we must have $\tau_{h,\delta_{h,j+1}} \geq p^{h+ts+s}$. Thus for any $0 \leq h \leq s-1$ and $1 \leq j \leq m-2$ we have

$$\sum_{i=0}^{\delta_{h,j}}\tau_{h,i} \leq \frac{1}{p^{s-1}}\tau_{h,\delta_{h,j+1}}. \tag{7.4}$$

Combining 7.3 and 7.4 we deduce

$$\sum_{j=1}^{m-2} \tilde{Z}_j - \tilde{O}_j \leq \frac{1}{p^{s-1}} \sum_{h=\alpha}^{\beta} \sum_{j=1}^{m-2} \tau_{h,\delta_{h,j+1}}$$

$$= \frac{1}{p^{s-1}} \sum_{h=\alpha}^{\beta} \sum_{j=2}^{m-1} \tau_{h,\delta_{h,j}}$$

By Lemma 6.2(e), the nonzero summands among $\tau_{h,\delta_{h,2}}, \ldots, \tau_{h,\delta_{h,m-2}}$ are all distinct powers of $p$. Thus the sum of these terms is $\leq \tau_{h,\delta_{h,m-1}}$ (with equality if and only if $\tau_{h,\delta_{h,m-1}} = 0$). This implies

$$\sum_{j=1}^{m-2} \tilde{Z}_j - \tilde{O}_j \leq \frac{2}{p^{s-1}} \sum_{h=\alpha}^{\beta} \tau_{h,\delta_{h,m-1}} \leq \sum_{h=\alpha}^{\beta} \tau_{h,\delta_{h,m-1}}$$

since we are assuming $s \geq 2$ and $p \geq 2$.

By Lemma 6.2(a), we have $\delta_{\beta,m-1} = u_\beta$. Thus $\tau_{\beta,\delta_{\beta,m-1}}$ is the largest term in $\tau_\beta(N)$. We have by Proposition 5.1(c) that $O_m = p^k$. Recall that $k \equiv \beta \pmod{s}$. Since $O$ is $\tau$-monotonic (Lemma 3.1) we have $\tau_{\beta,\delta_{\beta,m-1}} = p^k$. ∎

Combining Lemma 7.1 with (7.1) we deduce

$$wt(Z) - wt(O) \geq Z_m - \left(2p^k + \sum_{h=\alpha}^{\beta-1} \tau_{h,\delta_{h,m-1}}\right). \tag{7.5}$$

Now we give some lower bounds on $Z_m$.

**Lemma 7.2.** $\deg_p(Z_m) > k$.

**Proof.** Recall that $(G_1, \ldots, G_m)$ is the greedy element of $V_m(N)$. It is sufficient to show $\deg_p(Z_m) \geq \deg_p(G_m)$ since $\deg_p(G_m) > k$ (Proposition 5.1(b)). Define $\gamma$ by $\gamma \equiv \deg_p(G_m) \pmod{s}$, $0 \leq \gamma \leq s-1$. Then $u_\gamma > 0$ where $\bar{u} = [u_0, \ldots, u_{s-1}]^t = \Gamma(N)$. By Lemma 6.1(a) we have $\alpha \leq \gamma < \beta$. Since $u_\gamma > 0$, Parts (b) and (c) of Lemma 6.2 imply that $\delta_{\gamma,m-1}$ is strictly less than $u_\gamma$. Recall that

$$\Gamma(Z_m) = \bar{b}_m = [b_{0,m}, \ldots, b_{s-1,m}]^t = \bar{u} - \bar{\delta}_{m-1}.$$

Thus $b_{\gamma,m} > 0$. Since $b_{\gamma,m}$ is the length of $\tau_\gamma(Z_m)$ and since $Z$ is $\tau$-monotonic, the sequence $\tau_\gamma(Z_m)$ contains the largest $b_{\gamma,m}$ terms in $\tau_\gamma(N)$. Since $p^{\deg_p(G_m)}$ is a term in $\tau_\gamma(N)$, we have $\deg_p(Z_m) \geq \deg_p(G_m)$. ∎

**Lemma 7.3.** *If $\sum_{h=\alpha}^{\beta-1} \tau_{h,\delta_{h,m-1}}$ is positive then*

$$\deg_p(Z_m) > \deg_p\left(\sum_{h=\alpha}^{\beta-1} \tau_{h,\delta_{h,m-1}}\right) + 1.$$

**Proof.** Suppose $\sum_{h=\alpha}^{\beta-1} \tau_{h,\delta_{h,m-1}}$ is positive and let $r = \deg_p\left(\sum_{h=\alpha}^{\beta-1} \tau_{h,\delta_{h,m-1}}\right)$. Since the non-zero terms $\tau_{h,\delta_{h,m-1}}$ are all distinct powers of $p$ there is no carryover of $p$-adic digits in their sum. Thus $\tau_{\eta,\delta_{\eta,m-1}} = p^r$ for some $\alpha \leq \eta < \beta$ and we must have $\delta_{\eta,m-1} > 0$. Lemma 6.2(b) and (c) imply

$$b_{\eta,m} = u_\eta - \delta_{\eta,m-1} \geq p - 1$$

where $[u_1, \ldots, u_m]^t = \Gamma(N)$. Thus there are at least $p-1$ terms in the sequence $\tau_\eta(Z_m)$ (since its length is $b_{\eta,m}$). Since $Z$ is $\tau$-monotonic, all of these terms are greater than or equal to $p^r$. Since there can be a maximum of $p-1$ occurrences of the same term appearing in all of $\tau_\eta(N)$, the largest term in $\tau_\eta(Z_m)$ must be $\geq p^{r+s}$. The result follows since we are assuming $s \geq 2$. ∎

Finally, we complete the proof of Theorem 1.2.

**Proposition 7.4.** *If $s \geq 2$ then $\Psi$ is empty.*

**Proof.** Set $Q := p^k + \sum_{h=\alpha}^{\beta-1} \tau_{h,\delta_{h,m-1}}$. By (7.5) it is sufficient to show $Z_m > Q + p^k$. Recall that for each $h$ if $\tau_{h,\delta_{h,m-1}}$ is not zero then it equals $p^c$ for some $c \equiv h \pmod{s}$. Since $k \equiv \beta \pmod{s}$ the $p$-digits of $Q$ consist of zeros and ones. Thus, if $p \geq 3$ then

$$\deg_p(Q + p^k) = \max\left\{k, \deg_p\left(\sum_{h=\alpha}^{\beta-1} \tau_{h,\delta_{h,m-1}}\right)\right\}.$$

Lemmas 7.2 and 7.3 imply that if $p \geq 3$ then $Z_m > Q + p^k$. So we assume $p = 2$. We see that

$$\deg_p(Q + p^k) \leq \max\left\{k+1, \deg_p\left(\sum_{h=\alpha}^{\beta-1} \tau_{h,\delta_{h,m-1}}\right) + 1\right\}.$$

By Lemma 7.3 we may assume $\deg_p(Q+p^k) = k+1 = \deg_p(Z_m)$. Recall that $Z$ is $\tau$-monotonic and that $\Gamma(Z_m) = \bar{u} - \bar{\delta}_{m-1}$ where $\bar{u} = [u_1, \ldots, u_m]^t = \Gamma(N)$. Since $k+1 = \deg_p(Z_m)$ we must have $u_{\beta+1} - \delta_{\beta+1,m-1} > 0$. But then Lemma 6.2(a) implies $\beta = s-1$ and $k+1 \equiv \alpha = 0 \pmod{s}$. Since the length of the sequence $\tau_h(Z_m)$ is $u_h - \delta_{h,m-1}$, Lemma 6.2(b) implies that if $\alpha < h < \beta$ and $u_h > 0$ then $\tau_h(Z_m)$ has at least one term $\tau_{h,u_h}$. Lemma 6.2(c) implies that $\tau_\alpha(Z_m)$ has at least two terms $\tau_{\alpha,u_\alpha}$ and $\tau_{\alpha,u_\alpha-1}$. Note that $\tau_{\alpha,u_\alpha} = p^{k+1}$ since $k+1 = \deg_p(Z_m)$. To sum up the last three sentences, we have

$$Z_m \geq p^{k+1} + \tau_{\alpha,u_\alpha-1} + \sum_{h=\alpha+1}^{\beta-1} \tau_{h,u_h}.$$

Subtracting $Q + p^k$ we get

$$Z_m - (Q + p^k) \geq \tau_{\alpha,u_\alpha-1} - \tau_{\alpha,\delta_{\alpha,m-1}} + \sum_{h=\alpha+1}^{\beta-1} \left(\tau_{h,u_h} - \tau_{h,\delta_{h,m-1}}\right).$$

By Lemma 6.2(c), we have $u_\alpha - \delta_{\alpha,m-1} \geq 2$. Thus $\tau_{\alpha,u_\alpha-1} - \tau_{\alpha,\delta_{\alpha,m-1}} > 0$ since, as $p = 2$, any power of $p$ can appear as a $\tau_{i,j}$ at most once. Clearly we have $\sum_{h=\alpha+1}^{\beta-1} \tau_{h,u_h} - \tau_{h,\delta_{h,m-1}} \geq 0$. Thus $Z_m > Q + p^k$ as desired. ∎